\documentclass[12pt,a4paper,reqno]{amsart}
\usepackage{amsmath}
\usepackage{lscape}
\usepackage{graphics}

\begin{document}
\newcommand{\be}{\begin{equation}}
\newcommand{\ee}{\end{equation}} 
\newcommand{\linespace}{\vspace{\baselineskip}}
\newcommand{\upline}{\vspace{-\abovedisplayskip}\vspace{-\baselineskip}}
\renewcommand{\epsilon}{\varepsilon}
\newcommand{\half}{\tfrac{1}{2}}
\newcommand{\hlf}[1]{\tfrac{#1}{2}}
\newcommand{\third}{\tfrac{1}{3}}
\newcommand{\thrd}[1]{\tfrac{#1}{3}}
\newcommand{\quarter}{\tfrac{1}{4}}
\newcommand{\sixth}{\tfrac{1}{6}}
\newcommand{\qqquad}{\quad\quad\quad}
\newcommand{\qqqquad}{\quad\quad\quad\quad}
\renewcommand{\implies}{\quad\Longrightarrow\quad}
\renewcommand{\iff}{\quad\Longleftrightarrow\quad}
\renewcommand{\(}{\left(}
\renewcommand{\)}{\right)}
\newcommand{\1}{\mathbf{1}}
\newcommand{\A}{\mathcal{A}}
\newcommand{\mr}{\mathbb{R}} 
\newcommand{\R}{\mathbb{R}}
\newcommand{\Z}{\mathbb{Z}}
\newcommand{\mc}{\mathbb{C}} 
\newcommand{\C}{\mathbb{C}}
\newcommand{\mh}{\mathbb{H}} 
\renewcommand{\H}{\mathbb{H}}
\newcommand{\mo}{\mathbb{O}} 
\newcommand{\OO}{\mathbb{O}}
\newcommand{\mk}{\mathbb{K}} 
\newcommand{\K}{\mathbb{K}}
\newcommand{\D}{\mathbb{D}}
\newcommand{\mj}{\mathbb{J}} 
\newcommand{\J}{\mathbb{J}}
\newcommand{\F}{\mathbb{F}}
\newcommand{\pr}{^\prime}
\newcommand{\spr}{^{\,\prime}}
\newcommand{\mcs}{\widetilde{\mathbb{C} }} 
\newcommand{\sC}{\widetilde{\mathbb{C} }}
\newcommand{\mhs}{\widetilde{\mathbb{H} }} 
\newcommand{\sH}{\widetilde{\mathbb{H} }}
\newcommand{\mos}{\widetilde{\mathbb{O} }} 
\newcommand{\sO}{\widetilde{\mathbb{O} }}
\newcommand{\mks}{\widetilde{\mathbb{K} }} 
\newcommand{\sK}{\widetilde{\mathbb{K} }}
\newcommand{\sF}{\widetilde{\mathbb{F} }}
\newcommand{\is}{\tilde{i }} 
\newcommand{\itil}{\tilde{i}}
\newcommand{\js}{\tilde{j}} 
\newcommand{\jtil}{\tilde{j}}
\newcommand{\ktil}{\tilde{k}}
\newcommand{\ks}{\tilde{k}} 
\newcommand{\ls}{\tilde{l }} 
\newcommand{\ltil}{\tilde{l }}
\newcommand{\ils}{\tilde{il }} 
\newcommand{\jls}{\tilde{jl }} 
\newcommand{\kls}{\tilde{kl }} 
\newcommand{\es}{\tilde{e }} 
\newcommand{\etil}{\tilde{e }}
\newcommand{\sa}{\mathfrak{sa} } 
\newcommand{\so}{\mathfrak{so} } 
\renewcommand{\o}{\mathfrak{o} } 
\newcommand{\su}{\mathfrak{su} } 
\renewcommand{\u}{\mathfrak{u}}
\newcommand{\sq}{\mathfrak{sq} } 
\renewcommand{\sp}{\mathfrak{sp} } 
\renewcommand{\sl}{\mathfrak{sl}} 
\newcommand{\usp}{\mathfrak{usp}}
\newcommand{\gl}{\mathfrak{gl}} 
\newcommand{\g}{\mathfrak{g} }  
\newcommand{\h}{\mathfrak{h} } 
\newcommand{\n}{\mathfrak{n} } 
\renewcommand{\k}{\mathfrak{k} } 
\newcommand{\p}{\mathfrak{p} } 
\newcommand{\U}{\operatorname{U}} 
\newcommand{\SO}{\operatorname{SO} } 
\newcommand{\SU}{\operatorname{SU} } 
\renewcommand{\O}{\operatorname{O} }
\newcommand{\SP}{\operatorname{Sp} } 
\newcommand{\SL}{\operatorname{SL}} 
\newcommand{\GL}{\operatorname{GL}} 
\newcommand{\Der}{\operatorname{Der}} 
\newcommand{\Aut}{\operatorname{Aut}} 
\newcommand{\Str}{\operatorname{Str}} 
\newcommand{\Con}{\operatorname{Con}} 
\newcommand{\End}{\operatorname{End}} 
\newcommand{\Cl}{\operatorname{Cl}} 
\newcommand{\ds}{\operatorname{\dot{+}}}
\newcommand{\dsum}{\operatorname{\dot{+}}}
\newcommand{\osum}{\operatorname{\oplus}}
\newcommand{\ox}{\otimes}
\renewcommand{\aa}{\alpha}
\newcommand{\bb}{\beta}
\renewcommand{\gg}{\gamma}
\renewcommand{\d}{\delta }
\newcommand{\x}{\bar{x} }
\newcommand{\y}{\bar{y} }
\newcommand{\s}{\bar{s} }
\newcommand{\diag}{\operatorname{diag}}
\newcommand{\Span}{\operatorname{Span}}
\newcommand{\ad}{\operatorname{ad}}
\newcommand{\sgn}{\operatorname{sgn}}
\newcommand{\tr}{\operatorname{tr}}
\newcommand{\per}{\operatorname{per}} 
\newcommand{\Tri}{\operatorname{Tri}}
\renewcommand{\Im}{\operatorname{Im}}
\renewcommand{\Re}{\operatorname{Re}}
\newcommand{\id}{\operatorname{id}}
\newcommand{\<}{\langle}
\renewcommand{\>}{\rangle}
\newcommand{\I}{\leavevmode\hbox{\rm\small1\kern-3.8pt\normalsize1}}
\theoremstyle{plain} \newtheorem{prop}{Proposition}
\theoremstyle{plain} \newtheorem{theorem}{Theorem}[section]
\theoremstyle{plain} \newtheorem{lemma}{Lemma}[section]
\theoremstyle{definition} \newtheorem{defn}{Definition}
\numberwithin{equation}{section}
\hfuzz1pc 
\vfuzz1.2pt 

28 February 2002 \\
\title{Magic squares and matrix models of Lie algebras}
\author{C.\ H.\ Barton}
\curraddr{Department of Mathematics\\ University of York\\ Heslington \\
York \\ YO10 5DD}
\email[C.\ H.\ Barton]{Chris.Barton@bigfoot.com}
\author{A.\ Sudbery}
\email[A.\ Sudbery]{as2@york.ac.uk}
\begin{abstract}
This paper is concerned with the description of exceptional simple Lie algebras
as octonionic analogues of the classical matrix Lie algebras. We review
the Tits-Freudenthal construction of the magic square, which includes
the exceptional Lie algebras as the octonionic case of a construction in
terms of a Jordan algebra of hermitian $3\times 3$ matrices (Tits) or
various plane and other geometries (Freudenthal). We present alternative
constructions of the magic square which explain its symmetry,
and show explicitly how the use of split composition algebras 
leads to analogues of the matrix Lie algebras $\su(3)$, $\sl(3)$ 
and $\sp(6)$. We adapt the magic square construction to include 
analogues of $\su(2)$, $\sl(2)$ and $\sp(4)$ for all real division
algebras. 
\end{abstract}

\maketitle

\section*{Contents}

\contentsline {section}{\tocsection {}{1}{Introduction}}{1}
\contentsline {section}{\tocsection {}{2}{Algebras: Notation}}{3}
\contentsline {section}{\tocsection {}{3}{The Tits construction}}{8}
\contentsline {section}{\tocsection {}{4}{Symmetrical constructions of the $n=3$ magic square}}{12}
\contentsline {subsection}{\tocsubsection {}{\quad 4.1}{The triality algebra $\operatorname {Tri}\mathbb {K}$ and $\operatorname {Der}H_3(\mathbb {K})$}}{12}
\contentsline {subsection}{\tocsubsection {}{\quad 4.2}{The Vinberg construction}}{19}
\contentsline {subsection}{\tocsubsection {}{\quad 4.3}{The triality construction}}{21}
\contentsline {section}{\tocsection {}{5}{The rows of the magic square}}{22}
\contentsline {section}{\tocsection {}{6}{Magic squares of $n\times n$ matrices}}{26}
\contentsline {subsection}{\tocsubsection {}{\quad 6.1}{The Santander-Herranz construction}}{27}
\contentsline {section}{\tocsection {}{7}{Maximal compact subalgebras}}{28}
\contentsline {section}{\tocsection {}{8}{The $n = 2$ magic square}}{34}
\contentsline {section}{\tocsection {Appendix}{A}{Matrix identities}}{40}
\contentsline {section}{\tocsection {Appendix}{}{References}}{43}

\linespace

\section{Introduction} \label{sec:intro} Semisimple Lie groups and Lie
algebras are normally discussed in terms of their root systems, which
makes possible a unified treatment and emphasises the common features of
their underlying structures. However, some classical
investigations~\cite{Weyl} depend on particularly simple matrix
descriptions of Lie groups, which are only available for the classical
groups. This creates a distinction between the
classical Lie algebras and the  exceptional ones, which  is maintained in some
more recent work (e.g.~\cite{Howe,Molev,Nazarov}). This paper is
motivated by the desire to give matrix descriptions of the
exceptional Lie algebras, assimilating them to the classical ones, with
a view to extending results like the Capelli identities to the
exceptional cases.

It has long been known~\cite{JacobsonELA} that most exceptional Lie algebras
are related to the exceptional Jordan algebra of $3\times 3$ hermitian
matrices with entries from the octonions, $\mo$. Here we show that this
relation yields descriptions of certain real forms of the complex Lie
algebras $F_{4}$, $E_{6}$ and $E_{7}$ which can be interpreted as
octonionic versions of the Lie algebras of, respectively, antihermitian
$3\times 3$ matrices, traceless $3\times 3$ matrices and symplectic 
$6\times 6$ matrices. To be precise, we define for
each alternative algebra $\mk$ a Lie algebra $\sa (3,\mk )$ such that
$\sa (3,\mc )=\su(3)$ and $\sa(3,\mo)$ is the compact real form of
$F_{4}$; a Lie algebra $\sl(3,\mk)$ which is equal to $\sl (3,\mc)$ for $\mk=\mc$
and a non-compact real form of $E_{6}$ for $\mk=\mo$; and a Lie algebra
$\sp(6,\mk)$ such that  $\sp(6,\mc)$ is the set of $6\times 6$ complex
matrices $X$ satisfying $X^{\dag}J=-JX$ (where J is an antisymmetric
real $6\times 6$ matrix and $X^{\dag}$ denotes the hermitian conjugate
of $X$), and such that $\sp(6,\mo)$ is a non-compact real form of
$E_{7}$.

Our definitions can be adapted to $2\times 2$ matrices to yield Lie
algebras $\sa(2,\mk)$, $\sl (2,\mk)$ and $\sp(4,\mk)$ reducing to
$\su(2),\sl(2,\mc)$ and $\sp(4,\mc)$ when $\mk=\mc$. These Lie algebras
are isomorphic to various pseudo-orthogonal algebras.

The constructions in this paper are all related to Tits's magic square
of Lie algebras~\cite{Tits62}. This is a construction of a
Lie algebra $T(\K,\J)$ for any alternative algebra $\mk$ and Jordan 
algebra $\mj$. If $\K = \K_1$ and $\mj=H_{3}(\mk_{2})$ is the Jordan algebra of $3\times
3$ hermitian matrices over another alternative algebra $\mk_{2}$, the
Jordan product being the anticommutator, this yields a Lie algebra
$L_{3}(\mk_{1},\mk_{2})$ for any pair of alternative algebras. Taking
$\mk_{1}$ and $\mk_{2}$ to be real division algebras, we obtain a
$4\times 4$ square of compact Lie algebras which (magically) is
symmetric and contains the compact real forms of $F_{4}$, $E_{6}$, $E_{7}$
and $E_{8}$. We will show that if the division algebra $\mk_{1}$ is
replaced by its split form $\mks_{1}$, one obtains a non-symmetric
square of Lie algebras whose first three rows are the sets of matrix Lie
algebras described above: \begin{align} \label{eqn:iso2}
    L_{3}(\mk ,\mr )&=\sa (3,\mk) \notag \\
    L_{3}(\mk, \mcs )&=\sl (3,\mk) \\
    L_{3}(\mk ,\mhs )&=\sp (6,\mk). \notag
\end{align}
We will also describe magic squares of Lie algebras based on $2\times 2$
matrices, which have similar properties.

The organisation of the paper is as follows. In
Section~\ref{notation} we establish notation, recall the
definitions of various kinds of algebra, and introduce our generalised
definitions of the Lie algebras $\sa(n,\K)$, $\sl(n,\K)$ and
$\sp(2n,\K)$. In Section~\ref{sec:Tits} we present Tits's general
construction of $T(\K,\J)$, show that appropriate (split) choices of $\K$
yield the derivation, structure and conformal algebras of $\J$, and
describe the magic squares obtained when $\J$ is a Jordan algebra of
$3\times 3$ hermitian matrices. Section~\ref{symmetry} presents two
alternative, manifestly symmetric, constructions of the magic square;
one, due to Vinberg, describes the algebras in terms of matrices over
tensor products of division algebras, while the other is based on
Ramond's concept \cite{Ramond} of the 
triality algebra. In Section~\ref{rows} we develop the description of
the rows of the non-compact magic square as unitary, special linear and
symplectic Lie algebras, and briefly describe Freudenthal's geometrical
interpretation. In Section~\ref{nxn} we discuss the extension of these
results from $n=3$ to general $n$, which is only possible for the
associative division algebras $\K=\R$, $\C$ and $\H$, and in
Section~\ref{2mat} we discuss the case $n=2$, with the octonions again
included. In an appendix we prove extensions to alternative algebras of
various matrix identities that are needed throughout the paper.

\section{Algebras: Notation}
\label{notation}
When dealing with a Lie algebra, we will use the notation $\ds$ for the 
direct sum of vector spaces. As well as making formulae
easier to read, this enables us to reserve the use of $\oplus$ to denote the
direct sum of Lie algebras, i.e. $L = M\oplus N$ implies that $[M,N]=0$
in $L$. For direct sums of a number of copies of a vector space we will
use a multiple (rather than power) notation, writing $nV = V\ds
V\ds \cdots V$ ($n$ times). For real vector spaces with no Lie algebra
structure but with a pseudo-orthogonal structure (a preferred bilinear
form, not necessarily positive definite), we replace
$\ds$ by $\osum$ to denote the internal direct sum of orthogonal
subspaces, or the external direct sum of pseudo-orthogonal spaces in which
$V\osum W$ has the inherited inner product making $V$ and $W$ orthogonal
subspaces.

Let $\mk$ be an algebra over $\mr$ with a non-degenerate quadratic form
$x\mapsto \left| x\right| ^{2}$ and associated bilinear form $\< x,y\>$.
If the quadratic form satisfies
\be
\label{eqn:qform}
    |xy|^{2}=|x|^{2}|y|^{2}, \qquad
\forall x,y\in\mk ,
\end{equation}
then $\K$ is a \textit{composition algebra}. We consider $\mr$ to be
embedded in $\mk$ as the set of scalar multiples of the identity
element, and denote by $\mk ^{\prime}$ the subspace of $\mk$ orthogonal
to $\mr$, so that $\mk=\mr \ds \mk^{\prime }$;
we write $x=\Re x + \Im x$ with $\Re x \in \mr$ and $\Im x \in \mk
^{\prime}$. It can then be shown \cite{KantorSolod} that the conjugation which fixes each
element of $\mr$ and multiplies every element of $\mk ^\prime $ by $-1$,
denoted $x\mapsto \overline{x}$, satisfies
\linespace
\be
    \overline{xy}=\overline{y}\, \overline{x}
\end{equation}
and
\be
    x\overline{x}=\left| x\right| ^{2}.
\end{equation}
The inner product in $\K$ is given in terms of this conjugation as 
\[
 \< x, y\> = \Re (x\overline{y}) = \Re (\overline{x}y). 
\]

In a composition algebra $\K$ we will generally adopt the
typographical convention that lower-case letters from the end of the
roman alphabet ($\ldots, x, y, z$) denote general elements of $\K$,
while elements of $\K\pr$ are denoted by letters from the beginning of
the roman alphabet ($a,b,c,\ldots$).
We use the notations $[x,y]$ and $[x,y,z]$ for the commutator and associator
\begin{align*}
    [x,y] &= xy -yx,\\
    [x,y,z] &=(xy)z-x(yz).
\end{align*}
These change sign if any one of their arguments is conjugated:
\[ 
[\overline{x}, y] = - [x,y], \qqquad [\overline{x},y,z] = - [x,y,z].
\]
Any composition algebra $\mk$ satisfies the \textit{alternative law}, i.e. the
associator is an alternating function of $x,y$ and $z$ \cite{Schafer66}.

A division algebra is an algebra in which
\begin{equation*}
    xy=0 \implies  x=0 \; \text{ or }\;  y=0.
\end{equation*}
This is true in a composition algebra if the quadratic form $|x|^2$ is
positive definite. By Hurwitz's Theorem~\cite{Schafer66}, the only such positive definite composition algebras are  $\mr ,\mc
,\mh $ and $\mo $. These algebras are obtained by the
Cayley-Dickson process~\cite{Schafer66}; the same process with
different signs yields \textit{split} forms
of these algebras. These are so called because the familiar
equation
\begin{equation*}
    i^{2}+1=0
\end{equation*}
in $\C$, $\H$ or $\OO$ is replaced in the split algebras $\mcs$, $\mhs$ and $\mos$ by
\begin{equation*}
    i^{2}-1=(i+1)(i-1)=0
\end{equation*}
i.e.\ the equation can be \textit{split} for at least one of the imaginary
basis elements (specifically: for the one imaginary unit of $\sC$, for two
of the three imaginary units of $\sH$, and for four of the seven
imaginary units of $\sO$). Unlike the positive
definite algebras $\mr ,\mc ,\mh $ and $\mo $, the
split forms $\mcs, \mhs $ and $\mos$ are not division algebras.

In any algebra we denote by $L_x$ and $R_x$ the maps of left and right
multiplication by $x$:
\[ 
L_x(y) = xy, \qqquad R_x(y) = yx.
\]

Matrix notation: $\1$ denotes the identity matrix (of a size which will be
clear from the context), $X^\prime$ denotes the traceless part of the
$n\times n$ matrix $X$:
\[
 X^\prime = X - \frac{\tr X}{n}\1, 
\]
and $X^{\dag }$ denotes the hermitian conjugate of the
matrix $X$ with entries in $\mk$, defined in analogy to the complex case
by
\begin{equation*}
    (X^{\dag})_{ij}=\overline{x_{ji}}.
\end{equation*}

Our notation for Lie algebras is that of \cite{Sudbery84}.
We use $\su (s,t)$ for the Lie algebra of the unimodular pseudo-unitary group,
\begin{equation*}
    \su (s,t)=\{ X \in \mc ^{n\times n} :X^{\dag}G+GX=0,\,\tr X =0\}
\end{equation*}
where $G=\diag (1,\dots ,1,-1,\dots ,-1)$ with $s$ $+$ signs and $t$ $-$ signs;
$\sq (n)$ for the Lie algebra of antihermitian quaternionic matrices $X$,
\begin{equation*}
    \sq (n)=\{ X\in \mh ^{n\times n}: X^{\dag}=-X\} ;
\end{equation*}
and $\sp (2n,\mk )$ for the Lie algebra of the symplectic group of $2n\times
2n$ matrices with entries in $\K = \R$ or $\C$,
\begin{equation*}
     \sp (2n,\mk )=\{ X\in \mk ^{2n\times 2n}:X^{\dag }J+JX=0\}
\end{equation*}
where $J=\begin{pmatrix} 0 & I_{n} \\ -I_{n} & 0 \\ \end{pmatrix}$. For
a general composition algebra $\K$, however, this set of matrices is not
closed under commutation. We will denote it by
\[
 Q_{2n}(\K) = \{ X\in \mk ^{2n\times 2n}:X^{\dag }J+JX=0\} 
\]
and its traceless subspace by $Q\pr_{2n}(\K)$.
We will see that this can be extended to a Lie algebra $\sp(2n,\K)$.
We also have $\so (s,t)$, the Lie algebra of the pseudo-orthogonal
group $\SO (s,t)$, given by
\begin{equation*}
    \so (s,t)=\{ X\in \mr ^{n\times n}:X^{T}G+GX=0\}
\end{equation*}
where $G$ is defined as before. 

We will write $\O (V,q)$ for the
group of linear maps of the vector space $V$ preserving the
non-degenerate quadratic form $q$, $\SO (V,q)$ for its unimodular (or
special) subgroup, and $\o(V,q)$ or $\so (V,q)$ for their common Lie algebra. We
omit $q$ if it is understood from the context. Thus for any division
algebra we have the group $\SO (\mk )$ and the Lie algebra $\so (\mk )$.

A \textit{Jordan algebra} $\mj$ is defined to be a commutative algebra (over a
field which in this paper will always be $\R$) in which all products satisfy the Jordan identity
\be\label{Jordan}
    (xy)x^{2}=x(yx^{2}).
\end{equation}

Let $M_{n}(\mk )$ be the set of all $n\times n$ matrices with entries in
$\mk$, and let  $H_{n}(\mk )$ and $A_{n}(\mk )$ be the sets of all
hermitian and antihermitian matrices with entries in $\mk $
respectively. We denote by $H_{n}^{\prime}(\mk )$, $A_{n}^{\prime}(\mk
)$ and $M_{n}^{\prime}(\mk )$  the subspaces of traceless
matrices  of  $H_{n}(\mk )$, $A_{n}(\mk )$ and $M_{n}(\mk )$ respectively.
We thus have $M_{n}(\mk )=H_{n}(\mk )\ds A_{n}(\mk )$ and
$M_{n}^{\prime}(\mk)=H_{n}^{\prime }(\mk )\ds A_{n}^{\prime }(\mk )$. We will use the
fact that $H_{n}(\mk )$ is a Jordan algebra for $\mk = \mr ,\mc$ and $\mh$
for all $n$ and for
$\mk=\mo$ when $n=2,3$~\cite{Sudbery84}, with the Jordan product as the anticommutator
\begin{equation*}
    X\cdot Y=XY+YX.
\end{equation*}
This is a commutative but non-associative product.

We denote the associative algebra of linear endomorphisms of a vector
space $V$ by $\End V$. The \textit{derivation} algebra, $\Der \A$, of any algebra
$\A$ is the Lie algebra
\be
\label{eqn:Der}
    \Der \A = \{D\in\End \A \mid  D(xy)=D(x)y + xD(y),\; \forall x,y\in \A \}
\end{equation}
with bracket given by the commutator. The derivation algebras of the four
positive definite composition algebras are as follows:
\begin{align}\label{DerK}
    \Der \mr &= \Der \mc =0; \\
    \Der \mh &= C(\mh ^{\prime})=\{ C_{a} \mid a \in \mh ^{\prime}
\} \text{ where } C_{a}(q)=aq-qa\\
       &\cong\su (2) \cong\so(3);\\
\Der \mo  &\text{ is a compact exceptional Lie algebra of type } G_{2}.
\end{align}

In both alternative algebras and Jordan algebras there are constructions
of derivations from left and right multiplication maps. In an
alternative algebra $\K$
\be
D_{x,y} = [L_x, L_y] + [L_x, R_y] + [R_x, R_y] \label{D-def}
\end{equation}
is a derivation for any $x,y\in\K$, also given by
\be
D_{x,y}(z) = [[x,y],z] - 3[x,y,z]. \label{Daction}
\end{equation}
It satisfies the Jacobi-like identity (\cite{Schafer66}, p.78)
\be
D_{[x,y],z} + D_{[y,z],x} + D_{[z,x],y} = 0. \label{DJacobi}
\end{equation}

In Section \ref{symmetry}, following Ramond~\cite{Ramond}, we will extend the
derivation algebra to the \textit{triality} algebra, which consists of
triples of linear maps $A,B,C: \A \mapsto \A$ satisfying
\be\label{eqn:Tri}
    A(xy) = (Bx)y + x(Cy),\qquad\forall x,y\in \A.
\end{equation}

The \textit{structure algebra} $\Str \A$ of any algebra $\A$ is
defined to be the Lie algebra generated by left and right multiplication
maps $L_a$ and $R_{a}$ for $a\in \A$. For a Jordan algebra with identity
this can be shown to be~\cite{Schafer66}
\be
\label{eqn:Str}
    \Str  \mathbb{J} = \Der \mathbb{J} \ds L(\mathbb{J})
\end{equation}
where $L(\mathbb{J} )$ is the set of all $L_x$ with $x\in \mathbb{J}$.
The Jordan axiom \eqref{Jordan} implies that the commutator $[L_x,L_y]$ 
is a derivation of $\J$ for
all $x,y\in \J$; thus the Lie algebra structure of $\Str\J$ is defined by
the statements that $\Der\J$ is a Lie subalgebra and
\begin{align*}
[D, L_x] &= L_{Dx} \qquad (D\in \Der\J,\; x\in\J),\\
[L_x, L_y] &= L_x L_y - L_y L_x \in \Der \J \qquad (x,y\in\J).
\end{align*}
We denote by $\Str ^{\prime} \mathbb{J}$ the quotient of $\Str\J$ by the
subspace of multiples of $L_e$ where $e$ is the identity of $\J$. Both 
$\Str\J$ and $\Str ^{\prime}\J$ have an involutive
automorphism $T\mapsto T^*$ which leaves $\Der\J$ fixed and multiplies
each element of $L(\J)$ by $-1$.

We also require another Lie algebra
associated with a Jordan algebra with identity, namely the conformal algebra as
constructed by Kantor (1973) and Koecher (1967). This is the vector space
\be
\label{eqn:Con}
    \Con  \mathbb{J} = \Str  \mathbb{J} \ds 2\mathbb{J}
\end{equation}
with brackets
\begin{align}
[T,\,(x,y)] =&\, (Tx, T^*y),\label{Conbrac1}\\
[(x,0),\,(y,0)] = 0& = [(0,x),\,(0,y)],\label{Conbrac2}\\
[(x,0),\,(0,y)] =&\, \half L_{xy} + \half[L_x, L_y].\label{Conbrac3}
\end{align}
where * is the involution described above.

When $\J$ is a Jordan algebra of symmetric or hermitian matrices, the Lie
algebras $\Der\J$, $\Str\pr\J$ and $\Con\J$ can be identified with
matrix Lie algebras:
\begin{align}
\Der H_n(\R) &= A_n\pr(\R) = \so(n),\\
\Der H_n(\C) &= A_n\pr(\C) = \su(n),\\
\Str\pr H_n(\K) &= M_n\pr(\K) = \sl (\K) \qqqquad 
(\K = \R \text{ or }\C),\\
\Con H_n(\K) &= Q_{2n}\pr(\K) = \sp (2n, \K) \qquad 
(\K = \R \text{ or }\C).
\end{align}
We adopt these as definitions of the following series of Lie algebras for
any composition algebra:

\begin{defn} If $\K$ is a real composition algebra and $n$ is a natural
number such that $H_n(\K)$ is a Jordan algebra,
\begin{align}
\sa(n,\K) &= \Der H_n(\K), \label{defsa}\\
\sl(n,\K) &= \Str\pr H_n(\K), \label{defsl}\\
\sp(2n,\K) &= \Con H_n(\K) \label{defsp}
\end{align}
where $J=\begin{pmatrix} 0 & I_{n} \\ -I_{n} & 0 \\ \end{pmatrix}$.
\end{defn}

Thus $\sa(n,\R) = \so(n)$, $\sa(n,\C) = \su(n)$, and we will see in
Section~\ref{rows} that $\sa(n,\H) = \sq(n)$. We will also find matrix
descriptions of the other quaternionic Lie algebras as follows:
\begin{align} 
\sl(n,\H) &= \{X\in\H^{n\times n}: \Re(\tr X) = 0\},\\
\sp(2n,\H) &= \{X\in\H^{2n\times 2n}: X^\dag J + JX = 0,\; \Re(\tr X) =
0\}.
\end{align}
We note that the standard notations \cite{Helgason, Gilmore}
for the quaternionic Lie algebras are
\begin{align*}
    \sa(n,\mh) &= \usp(2n)\\
    \sl(n,\mh) &= \su^*(2n)\\
    \sp(2n,\mh) &= \so^*(4n).
\end{align*}

\section{The Tits Construction}
\label{sec:Tits}

Let $\K$ be a real composition algebra and $(\J,\cdot)$ a real Jordan algebra
with identity $E$, and suppose $\J$ has an inner product $\<,\>$ satisfying 
\be
\<X,\, Y\cdot Z\> = \<X\cdot Y,\, Z\>. \label{Jassoc} 
\end{equation}
Let $\J\pr$ and $\K\pr$ be the subspaces of $\J$ and $\K$ orthogonal to
the identity, and let * denote the product on $\J\pr$ obtained from the
Jordan product by projecting back into $\J\pr$:
\[ 
A*B = A\cdot B - \frac{4}{n}\<A,B\>\1 \qquad \text{where } \frac{n}{4}=\<E,E\>
\]
(the notation is chosen to fit the case $\J = H_n(\K)$ with $X\cdot Y =
XY+YX$ and $\<X, Y\> = \half\tr (X\cdot Y)$; then $E = \half\1$ and $\<E,E\>=n/4$).
Tits defined a Lie algebra structure on the vector space
\be
T(\K,\J ) = \Der\K \ds \Der\J \ds \K\pr\ox\J\pr \label{Titsspace}
\end{equation}
with the usual brackets in the Lie subalgebra $\Der\K\oplus\Der\J$,
brackets between this and $\K\pr\ox\J\pr$ defined by the usual action of
$\Der\K\oplus\Der\J$ on $\K\pr\ox\J\pr$,
and further brackets
\be
[a\ox A, b\ox B ] = \frac{1}{n}\<A,B\>D_{a,b} 
- \<a,b\>[L_A,L_B] + \half[a,b]\ox (A*B) \label{Titsbrac} 
\end{equation}
where $a,b\in\K\pr$; $A,B\in\J\pr$; the square brackets on the
right-hand side denote commutators in $\K\pr$ and $\End\J$; and $D_{a,b}$
is the derivation of $\K\pr$ defined in (\ref{D-def}). For future
reference we sketch the proof of a slightly generalised version of Tits's theorem
(\cite{Tits62}; see also \cite{Freudenthal65}, \cite{Schafer66}).

\begin{theorem}\label{thm:Tits} \emph{(Tits)} The brackets \emph{(\ref{Titsbrac})} define a Lie
algebra structure on $T(\K,\J)$ if either $\K$ is associative or in $\J$ 
there is a cubic identity
\be
\frac{n}{6}X*(X\cdot X) = \<X,X\>X, \qquad \text{\emph{all} }X\in\J\pr. 
\label{cubic} 
\end{equation}
\end{theorem}

\begin{proof} The identity (\ref{Jassoc}) guarantees that derivations of
$\J$ are antisymmetric with respect to the inner product $\< ,\>$; since
the inner product in $\K$ is constructed from the multiplication,
the same applies in $\K$. It follows that the brackets
(\ref{Titsbrac}) are equivariant under the action of $\Der\K\oplus\Der\J$,
so that all Jacobi identities involving these derivations are satisfied.
Thus we need only consider the Jacobi identity between three elements of
$\K\pr\ox\J\pr$, namely the vanishing of 
\[ 
[[a\ox A, b\ox B],c\ox C] + [[b\ox B, c\ox C],a\ox A] + 
[[c\ox C, a\ox A], b\ox B].
\]

The component of this in $\Der\K$ is
\[ 
\frac{1}{8n}\<A*B,C\>\(D_{[a,b],c} + D_{[b,c],a} + D_{[c,a],b}\)
\]
which vanishes by (\ref{DJacobi}). The component in $\Der\J$ is
\[ 
-\half\<[a,b],c\>\big([L_{A\cdot B},L_C] + [L_{B\cdot C},L_A] + 
[L_{C\cdot A},L_B]\big)
\]
which vanishes by the polarisation of the Jordan axiom 
$[L_{X\cdot X}, L_X] = 0$ (obtained by putting $X = \lambda A + \mu B + \nu C$ and
equating coefficients of $\lambda\mu\nu$).

Finally, the component in $\K\pr\ox\J\pr$ is
\begin{align*}
Q = -&\frac{3}{n}[a,b,c]\ox\big(\<A,B\>C + \<B,C\>A + \<C,A\>B\big)\\
&+\(\<a,c\>b - \<a,b\>c + \quarter [[b,c],a]\)\ox A*(B\cdot C)\\
&+\(\<b,a\>c - \<b,c\>a + \quarter [[c,a],b]\)\ox B*(C\cdot A)\\
&+\(\<c,b\>a - \<c,a\>b + \quarter [[a,b],c]\)\ox C*(A\cdot B).
\end{align*}
Now $\<a,b\> = \Re(a\overline{b}) = -\half(ab + ba)$ since $b\in\K\pr$;
hence
\begin{align*}
4\big(\<a,c\>b - \<a,b\>c\big) &= - (ac+ca)b - b(ac+ca) + (ab+ba)c + c(ab+ba)\\
&= -[[b,c],a] + 2[a,b,c]
\end{align*}
and so 
\begin{multline*}
Q = [a,b,c]\ox \left\{-\frac{3}{n}\big(\<A,B\>C + \<B,C\>A + \<C,A\>B\big)\right.\\
\left.+\frac{1}{2}\big(A*(B\cdot C) + B*(C\cdot A) + C*(A\cdot B)\big)\right\}.
\end{multline*}
If $\K$ is associative, the first factor vanishes; if the identity
(\ref{cubic}) holds in $\J$, then polarising it shows that the second
factor vanishes.
\end{proof}

Taking $\K$ to be $\R$ or one of the split composition algebras $\sC$ or
$\sH$ gives three of the Lie algebras associated with $\J$ defined in
Section \ref{notation}:

\begin{theorem}\label{Jmatrix} For any Jordan algebra $\J$,
\begin{align}T(\R,\J) &\cong \Der\J,\\
T(\sC, \J) &\cong \Str\pr\J,\\
T(\sH,\J) &\cong \Con\J.
\end{align}
\end{theorem}

\begin{proof} Since $\Der\R = \R\pr = 0$, the first statement is true by
definition. 

For the second statement, we define $\theta : T(\sC,\J)\rightarrow\Str\pr\J$
to be the identity on $\Der\J$, and on $\sC\pr\ox\J\pr$
\[
\theta(\itil\ox A) = L_A.
\]
Since $\Der\sC = 0$ and $\sC\pr$ is spanned by $\itil$ which satisfies
$\<\itil,\itil\> = -1$, this is an isomorphism between $T(\sC,\J)$ and
the subspace of $\Str\J$ spanned by $\Der\J$ and the multiplication maps
$L_A$ with $A\in\J$, i.e.\ the subspace
$\Der\J\dsum L(\J\pr)\cong\Str\pr\J$ by (\ref{eqn:Str}).

In the third statement we have 
\begin{align*}
\Con\J &\cong \Str\J\dsum 2\J\\
&\cong \Der\J \dsum 3\J
\end{align*}
which is isomorphic to $T(\sH,\J)$ as a vector space since
\begin{align*}
T(\sH,\J) &= \Der\J \dsum \sH\pr\ox\J\pr \dsum \Der\sH \\
&\cong \Der\J \dsum 3\J\pr \dsum C(\sH\pr)\\
&\cong \Der\J \dsum 3\J
\end{align*}
since $C(\sH\pr)$ is 3-dimensional. Taking
the multiplication in $\sH$ to be given by 
\begin{align*}
\itil ^2 = \jtil ^2 &= 1, \quad \quad \ktil^2 = -1;\\
\itil \jtil = -\jtil\itil = -\ktil, \quad \jtil\ktil &= -\ktil\jtil = \itil,
\quad \ktil \itil = - \itil\ktil = \jtil,
\end{align*}
we define $\phi : T(\sH, \J)\rightarrow \Con\J$ by
\begin{align*}
\phi(D) &= D \qquad (D\in\Der\J),\\
\phi(C_{\itil}) &= 2L_1 \in \Str\J,\\
\phi(C_{\jtil}) &= 2(1,1) \in 2\J,\\
\phi(C_{\ktil}) &= -2(1,-1) \in 2\J,\\
\phi(\itil\ox A) &= L_A \in\Str\J,\\
\phi(\jtil\ox A) &= (A,A) \in 2\J,\\
\phi(\ktil\ox A) &= (-A,A) \in 2\J \quad (A\in\J\pr).
\end{align*}
It is straightforward to check that this is a Lie algebra isomorphism.
\end{proof}

Tits obtained the magic square of Lie algebras by taking the Jordan
algebra $\J$ to be the algebra of hermitian $3\times 3$ matrices over a
second composition algebra, defining
\[
L_3(\K_1,\K_2) = T(\K_1,H_3(\K_2)).
\]
The inner product in $H_3(\K_2)$ is given by $\<X,Y\> = \half\tr(X\cdot Y)$.
This yields the Lie algebras whose complexifications are
\begin{center}
\vspace{0.5cm} \begin{tabular}{|c c||c|c|c|c|}
\hline
 &$\K_2$  & $\mr$ & $\mc$  & $\mh$ & $\mo$ \\
$\K_1$&&&&&\\
 \hline \hline
  $\mr$ && $A_1$ & $A_2$ & $C_3$ & $F_4$ \\
  \hline
   $\mc$ && $A_2$ & $A_2 \oplus A_2$ & $A_5$ & $E_6$ \\
   \hline
    $\mh$ && $C_3$ & $A_5$ & $D_6$ & $E_7$ \\
    \hline
    $\mo$ && $F_4$ & $E_6$ & $E_7$ & $E_8$ \\
    \hline
    \end{tabular} \vspace{0.5cm}
    \end{center}
i.e.\ the Lie algebras with compact  real forms
\begin{center}
\vspace{0.5cm} \begin{tabular}{|c||c|c|c|c|}
\hline
  & $\mr$ & $\mc$  & $\mh$ & $\mo$ \\
 \hline \hline
  $\mr$ & $\so(3)$ & $\su(3)$ & $\sq(3)$ & $F_4$ \\
  \hline
   $\mc$ & $\su(3)$ & $\su(3) \oplus \su(3)$ & $\su(6)$ & $E_6$ \\
   \hline
    $\mh$ & $\sq(3)$ & $\su(6)$ & $\so(12)$ & $E_7$ \\
    \hline
    $\mo$ & $F_4$ & $E_6$ & $E_7$ & $E_8$ \\
    \hline
    \end{tabular}. \vspace{0.5cm}
    \end{center}
The striking properties of this square are (a) its symmetry and (b) the
fact that four of the five exceptional Lie algebras occur in its last
row. The fifth exceptional Lie algebra, $G_{2}$, can be included by adding an extra row
corresponding to the Jordan algebra
$\mr$. The explanation of the symmetry property is the subject of
the following section.

If one of the composition algebras is split, the magic square
$L_3(\sK_1,\K_2)$, according to Theorem \ref{Jmatrix}, contains matrix Lie
algebras as follows:
\begin{center}\vspace{0.5cm}
{\renewcommand{\arraystretch}{1.4}
\begin{tabular}{|c||c|c|c|c|}
\hline
     & $\mr$ & $\mc$  & $\mh$ & $\mo$ \\
 \hline \hline
  $\Der H_3(\mk )\cong L_3(\mk ,\mr)\cong \su (3,\mk)$ & $\so(3)$ & $\su(3)$ & $\sq(3)$ & $F_4(52)$ \\
  \hline
   $\Str ^\prime  H_3(\mk )\cong L_3(\mk ,\mcs ) \cong \sl (3,\mk )$ & $\sl(3,\mr )$ &
$\sl(3,\mc )$ & $\sl(3,\mh )$ & $E_6(26)$ \\
   \hline
    $\Con  H_3(\mk )\cong L_3(\mk ,\mhs ) \cong \sp (6,\mk )$ & $\sp(6,\mr )$ & $\su(3,3)$ &
$\sp(6,\mh )$ & $E_7(25)$ \\
    \hline\vspace{2pt}
    $L_3(\mk ,\mos) $ & $F_4(-4)$ & $E_6(-2)$ & $E_7(5)$ & $E_8(24)$ \\
    \hline
    \end{tabular} \vspace{0.5cm}}
    \end{center}

\noindent where the real forms of the exceptional Lie algebras in the last row and
column are labelled by the signatures of their Killing forms. These can be
also be identified by their maximal compact subalgebras
as follows:
\begin{center}
\label{maximal}
\vspace{0.5cm} \begin{tabular}{|c|c|}
\hline
    Exceptional Lie Algebra & Maximal Compact Subalgebra \\
\hline
    $F_4(52)$ & $F_4$\\
\hline
    $E_6(26)$ & $F_4$ \\
\hline
    $E_7(25)$ & $E_{6} \oplus \so (2)$ \\
\hline
    $E_8(24)$ & $E_{7} \oplus \so (3)$ \\
\hline
    $F_4(-4)$ & $\sq (3) \oplus \so (3)$ \\
\hline
    $E_6(-2)$ & $\su (6) \oplus \so (3)$ \\
\hline
    $E_7(5)$ & $\so (12) \oplus \so (3)$ \\
\hline
    $E_8(24)$ & $E_7\oplus \so(3)$\\
\hline
\end{tabular} \vspace{0.5cm}.
\end{center}

In section 5 we will give a general explanation of the close relationship between the
maximal compact subalgebras of the algebras in one line of the split
magic square
$L_3(\K_1,\sK_2)$ and the algebras in the preceding line of the
compact magic square $L_3(\K_1,\K_2)$. We will use the same method to
identify the maximal compact subalgebras of the doubly split magic
square $L_3(\sK_1,\sK_2)$.

\section{Symmetrical Constructions of the $n=3$ Magic Square}
\label{symmetry}

In this section we present two alternative constructions of Tits's
magic square which are manifestly symmetric between the two composition
algebras $\K_1, \K_2$. These are the constructions of Vinberg
\cite{Onishchik} and a construction using the triality algebra based on a
suggestion of Ramond \cite{Ramond}. We begin by exploring the structure
of the Lie algebras associated with the Jordan algebra $H_3(\K)$. 

\subsection{The triality algebra $\Tri\K$ and $\Der H_3(\K)$}
\label{Trider}

\begin{defn}
Let $\mk $ be a composition algebra over $\mr$. The \emph{triality
algebra} of $\mk $ is defined to be
\begin{equation}
\label{eqn:triality}
    \Tri \mk = \{ (A,B,C) \in 3\so (\mk ): A(xy)=x(By)+(Cx)y,\;
\forall x,y\in \mk \}.
\end{equation}
\end{defn}

The structure of the triality algebras $\Tri\mk$ can be analysed in a
unified way as follows:

\begin{lemma}\label{TriK}
For any composition algebra $\mk$,
\begin{equation*}
    \Tri \mk = \Der \mk \ds 2\mk ^{\prime}
\end{equation*}
in which $\Der \mk $ is a Lie subalgebra and the other brackets are
\begin{align*}
    [D,(a,b)] &= (Da, Db) \in 2\mk ^{\prime} \\
    [(a,0),(b,0)] &= \tfrac{2}{3}D_{a,b}+
\left( \tfrac{1}{3}[a,b],-\tfrac{2}{3}[a,b]\right), \\
    [(a,0),(0,b)] &= \tfrac{1}{3}D_{a,b}-\left( \tfrac{1}{3}[a,b],\tfrac{1}{3}[a,b]\right),
\\
    [(0,a),(0,b)] &= \tfrac{2}{3}D_{a,b}+\left(
-\tfrac{2}{3}[a,b],\tfrac{1}{3}[a,b]\right).
\end{align*}
\begin{proof}
Define $T:\Der \mk \ds 2\mk\rightarrow \Tri \mk $ by
\begin{multline}
\label{eqn:Teq}
    T(D,a,b)=
    (D+L_{a}-R_{b},\, D-L_{a+b}-R_{b},\, D+L_{a}+R_{a+b}).
\end{multline}
The alternative law guarantees that the right-hand side belongs to
$\Tri\K$; the Lie algebra isomorphism property follows from the brackets
\begin{align*}
    [L_{x},L_{y}]&=
\tfrac{2}{3}D_{x,y}+\tfrac{1}{3}L_{[x,y]}+\tfrac{2}{3}R_{[x,y]} \\
    [L_{x},R_{y}]&=
-\tfrac{1}{3}D_{x,y}+\tfrac{1}{3}L_{[x,y]}-\tfrac{1}{3}R_{[x,y]} \\
    [R_{x},R_{y}]&=
\tfrac{2}{3}D_{x,y}-\tfrac{2}{3}L_{[x,y]}-\tfrac{1}{3}R_{[x,y]}
\end{align*}
and the expression of the inverse map $T^{-1}:\Tri \mk \rightarrow \Der \mk
\ds 2\mk$ as
\begin{equation*}
    T^{-1}(A,B,C)=(A-L_{a}+R_{b},a,b)
\end{equation*}
where $a=\tfrac{1}{3}B(1)+\tfrac{2}{3}C(1)$ and
$b=-\tfrac{2}{3}B(1)-\tfrac{1}{3}C(1)$.
\end{proof}
\end{lemma}

Note that $\Tri\K$ has a subalgebra $\Der\K\dsum\K\pr$ in which the
second summand is the diagonal subspace of $2\K\pr$, containing the
elements $(a,a)$. Identifying $(a,a)\in2\K\pr$ with $a\in\K\pr$, the
brackets between elements of $\K\pr$ in this subalgebra are given by 
\[
 [a,b] = 2D_{a,b} - [a,b].
\]
The derivation and triality algebras of the four real division algebras,
together with this intermediate algebra, are tabulated below.

\begin{center}
\vspace{0.5cm} \begin{tabular}{|c||c|c|c|}
\hline
 $\K$ & $\Der\K$ & & $\Tri\K$ \\
 \hline \hline
  $\R$ & 0 & 0 & 0 \\
  \hline
   $\C$ & 0 & $\so(2)$ & $\so(2)\oplus\so(2)$ \\
   \hline
    $\H$ & $\so(3)$ & $\so(3)\oplus\so(3)$ & $\so(3)\oplus\so(3)\oplus\so(3)$ \\
    \hline
    $\OO$ & $G_2$ & $\so(7)$ & $\so(8)$ \\
    \hline
    \end{tabular} \vspace{0.5cm}
    \end{center}

The identification $\Tri\OO = \so(8)$ is a form of the principle of
triality \cite{Porteous}.

We will now show that any two elements $x,y\in\mk$ have an element of
Tri$\mk$ assicated with them, of which the first component is
the generator of rotations in the plane of $x$ and $y$, defined as
\begin{equation}\label{S-def}
S_{x,y}(z) = \< x,z\> y - \< y,z\> x.
\end{equation}
\begin{lemma}
\label{lemma:TTri}
For any $x,y\in \mk$, let
\begin{equation*}
    T_{x,y}=(4S_{x,y},\,
R_{y}R_{\bar{x}}-R_{x}R_{\bar{y}},\, L_{y}L_{\bar{x}}-L_{x}L_{\bar{y}}).
\end{equation*}
Then $T_{x,y}\in \Tri \mk $.

\begin{proof}
Write the action of $S_{x,y}$ as
\begin{align}
    2S_{x,y}(z) &= (x\bar{z}+z\bar{x})y-x(\bar{z}y+\bar{y}z) \label{eqn:cat}\\
    &= -[x,y,z]+z(\bar{x}y)-(x\bar{y})z \label{eqn:mouse}
\end{align}
using the alternative law and the relation $[x,y,\bar{z}]=-[x,y,z]$.
Since $\Re (\bar{x}y)=\Re (x\bar{y})$, we can write the last two terms as
\begin{align}
    z(\bar{x}y)-(x\bar{y})z &= z\Im (\bar{x}y)-\Im (x\bar{y})z \label{eqn:dog} \\
    &=\tfrac{1}{2}z(\bar{x}y-\bar{y}x)-\tfrac{1}{2}(x\bar{y}-y\bar{x})z.
 \label{eqn:budgie}
\end{align}
Now, by equation~\eqref{D-def}, we have
\[
\label{eqn:Sdef2}
    S_{x,y}=\tfrac{1}{6}D_{x,y}+L_{a}-R_{b}
\]
with
\begin{align*}
    a &= -\tfrac{1}{6}[x,y]-\tfrac{1}{4}(x\bar{y}-y\bar{x}) \, \in \mk
^{\prime} \\
    b &= -\tfrac{1}{6}[x,y]-\tfrac{1}{4}(\bar{x}y-\bar{y}x) \, \in \mk
^{\prime }.
\end{align*}
Hence by Lemma\ref{TriK} there is an element $(A,B,C) \in \Tri \mk $
with $A=S_{x,y}$ and
\begin{align*}
    B &= \tfrac{1}{6}D_{x,y}-L_{a+b}-R_{b}= S_{x,y}-L_{2a+b}, \\
    C &= \tfrac{1}{6}D_{x,y}+L_{a}+R_{a+b}= S_{x,y}+R_{a+2b}.
\end{align*}
Writing $[x,y]=-\tfrac{1}{2}([\bar{x},y]+[x,\bar{y}])$ gives
\begin{align*}
    a+2b&=\tfrac{1}{4}(\bar{y}x-\bar{x}y), \\
    2a+b &=\tfrac{1}{4}(y\bar{x}-x\bar{y});
\end{align*}
thus equations~\eqref{eqn:cat} and~\eqref{eqn:dog} imply that
\[
    S_{x,y} = -\tfrac{1}{2}A_{x,y}-R_{a+2b}+L_{2a+b}
\]
where $A_{x,y}(z) = [x,y,z]$. Hence
\begin{align*}
    Cz &= -\tfrac{1}{2}[x,y,z]+\tfrac{1}{4}(y\bar{x}-x\bar{y})z \\
    &= -\quarter [y,\overline{x},z] + \quarter[x,\overline{y},z] +
\quarter (y\overline{x} - x\overline{y})z\\
    &=\tfrac{1}{4}y(\bar{x}z)-\tfrac{1}{2}x(\bar{y}z),
\end{align*}
i.e.
\begin{equation*}
    C=\tfrac{1}{4}(L_{y}L_{\bar{x}}-L_{x}L_{\bar{y}}).
\end{equation*}
Similarly,
\begin{equation*}
    B=\tfrac{1}{4}(R_{y}R_{\bar{x}}-R_{x}R_{\bar{y}}).
\end{equation*}
Thus $(4S_{x,y},4C,4B)=T_{x,y}$, which is therefore an element of $\Tri \mk $.
\end{proof}
\end{lemma}

Define an automorphism of $\Tri \mk $ as follows. For any linear map
$A:\mk \rightarrow \mk$ let $\overline{A}=KAK$ where $K:\mk \rightarrow \mk$ is
the conjugation $x\mapsto \bar{x}$, i.e.
\begin{equation*}
    \overline{A}(x)=\overline{A(\bar{x})}.
\end{equation*}
\begin{lemma}
\label{lemma:splodge}
Given $T=(A,B,C) \in \Tri \mk$, let
\begin{equation*}
    \theta (T) = (\overline{B},C,\overline{A}).
\end{equation*}
Then $\theta (T) \in \Tri \mk$ and $\theta$ is a Lie algebra
automorphism.
\begin{proof}
By Lemma~\ref{TriK}, $T=T(D,a,b)$ for some $D\in \Der \mk$ and $a,b
\in \mk ^{\prime}$. Then
\begin{align*}
    A&= D+L_{a}-R_{b} \\
    B&= D-L_{a+b}-R_{b} \\
    C&= D+L_{a}+R_{a+b}.
\end{align*}
It follows that
\begin{equation*}
      \overline{B}=D+R_{a+b}+L_{b}=D+L_{a^{\prime}}-R_{b^{\prime}}
\end{equation*}
with $a^\prime = b$, $b^\prime = -a-b$. This is the first component of
the triality $T^{\prime}=(A^{\prime},B^{\prime},C^{\prime})$, where

\begin{align*}
B^{\prime} &=D-L_{a^{\prime} + b^{\prime}}-R_{b^{\prime}} = C,\\
C^{\prime} &= D+L_{a^{\prime}}+R_{a^{\prime} + b^{\prime}} = \overline{A}
\end{align*}
i.e.\ $T^{\prime}=(\overline{B},C,\overline{A})=\theta(T)$. It is clear
that $\theta$ is a Lie algebra automorphism.
\end{proof}
\end{lemma}

\begin{theorem}
\label{theorem:banana}
For any composition algebra $\mk$,
\begin{equation*}
    \Der H_{3}(\mk)=\Tri \mk \ds 3\mk
\end{equation*}
in which $\Tri \mk $ is a Lie subalgebra; the brackets in $[ \Tri
\mk ,3\mk ]$ are
\begin{equation}
\label{eqn:1}
    [T,F_{i}(x)]=F_{i}(T_{i}x)\in 3\mk ,
\end{equation}
if $T=(T_{1},\overline{T}_{2},\overline{T}_{3}) \in \Tri \mk$ and
$F_{1}(x)+F_{2}(y)+F_{3}(z)=(x,y,z) \in 3\mk $; and the brackets in
$[\Tri \mk ,\Tri \mk ]$ are given by
\begin{align}
\label{eqn:3}
    [F_{i}(x),F_{i}(y)]&= \theta ^{1-i}(T_{x,y}) \in \Tri \mk ,\\
\label{eqn:2}
    [F_{i}(x),F_{j}(y)]&=F_{k}(\bar{y}\bar{x}) \in 3\mk,
\end{align}
if $x,y \in \mk$ and $(i,j,k)$ is a cyclic permutation of $(1,2,3)$.
\begin{proof}
Define elements $e_{i},P_{i}(x)$ of $H_{3}(\mk)$ (where $i=1,2,3;\,  x \in
\mk $) by
\begin{equation}
\label{eqn:Jdef}
\begin{pmatrix} \alpha & z & \bar{y} \\ \bar{z} & \beta & x \\ y &
\bar{x} & \gamma \end{pmatrix} = \alpha e_{1}+\beta e_{2}+\gamma e_{3} +
P_{1}(x)+P_{2}(y)+P_{3}(z)
\end{equation}
for $\alpha , \beta , \gamma \in \mr$ and $x,y,z \in \mk $. The Jordan
product in $H_{3}(\mk)$ is given by
\begin{subequations}
\begin{align}
    e_{i}\cdot e_{j}&=2\delta _{ij} e_{i}  \label{eqn:J1} \\
    e_{i}\cdot P_{j}(x) &= (1-\delta _{ij})P_{j}(x) \label{eqn:J2}\\
    P_{i}(x)\cdot P_{i}(y) &= 2\<x,y\>(e_{j}+e_{k}) \label{eqn:J3} \\
    P_{i}(x)\cdot P_{j}(y) &= P_{k}(\bar{y}\, \bar{x})
\label{eqn:J4}
\end{align}
\end{subequations}
where in each of the last two equations $(i,j,k)$ is a cyclic
permutation of $(1,2,3)$.

Now let $D:H_{3}(\mk)\to H_{3}(\mk)$ be a derivation of this algebra.
First suppose that
\begin{equation*}
    De_{i}=0, \quad i=1,2,3.
\end{equation*}
Then
\begin{align*}
    e_{i}\cdot DP_{i}(x) &=0, \\
    e_{i}\cdot DP_{j}(x) &= DP_{j}(x) \quad \text{if } i \neq j.
\end{align*}
Thus $DP_{j}(x)$ is an eigenvector of each of the multiplication
operators $L_{e_{i}}$, with eigenvalue $0$ if $i=j$ and $1$ if $i\neq
j$. It follows that
\begin{equation}
\label{eqn:derp}
    DP_{j}(x)=P_{j}(T_{j}x)
\end{equation}
for some $T_{j}:\mk \to \mk $. Now
\begin{equation*}
    DP_{j}(x)\cdot P_{j}(y) +P_{j}(x)\cdot DP_{j}(y)=0
\end{equation*}
gives $T_{j} \in \so (\mk )$; and the derivation property of $D$ applied
to \eqref{eqn:J4} gives
\begin{equation*}
    T_{k}(\bar{y}\bar{x})=\bar{y}(\overline{T_{i}x})+(\overline{T_{j}y})\bar{x}
\end{equation*}
i.e.\ $(T_{k},\overline{T_{i}},\overline{T_{j}}) \in \Tri \mk $ and
therefore $(T_{1},\overline{T_{2}},\overline{T_{3}})\in \Tri \mk $.

If $De_{i}\neq 0$, then from equation~\eqref{eqn:J1} with $i=j$,
\begin{equation*}
2e_{i}\cdot De_{i}=2De_{i}
\end{equation*}
so $De_{i}$ is an eigenvector of the multiplication $L_{e_{i}}$ with
eigenvalue $1$, i.e.\ $De_{i} \in P_{j}(\mk )+P_{k}(\mk )$ where $(i,j,k)$
are distinct. Write
\begin{equation*}
    De_{i} = P_{j}(x_{ij} )+P_{k}(x_{ik} );
\end{equation*}
then equation~\eqref{eqn:J1} with $i\neq j$ gives
\begin{equation*}
    e_{i}\cdot P_{k}(x_{jk})+e_{i}\cdot
P_{i}(x_{ji})+P_{j}(x_{ij})\cdot e_{j}+P_{k}(x_{ik})\cdot e_{j} =0.
\end{equation*}
Hence
\begin{equation*}
    P_{k}(x_{jk}+x_{ik})=0.
\end{equation*}
It follows that the action of any derivation on the $e_{i}$ must be of
the form $De_i = F_{1}(x)+F_{2}(y)+F_{3}(z)$ where
\begin{align}
\label{eqn:Fone}
    F_{i}(x)e_{i} &=0 \notag \\
    F_{i}(x)e_{j}=-F_{i}(x)e_{k}&=P_{i}(x),
\end{align}
$(i,j,k)$ being a cyclic permutation of $(1,2,3)$. Hence 
\[
\Der H_{3}(\mk)\subseteq \Tri \mk \oplus 3\mk.
\]
To show that such derivations $F_{i}(x)$ exist and therefore that the
inclusion just obtained is an equality, consider the
operation of commutation with the matrix
\begin{align*}
    X&=\begin{pmatrix} 0 & -z & \bar{y} \\ \bar{z} & 0 & -x \\ -y &
\bar{x} & 0 \end{pmatrix} \\
&\\
    &= X_{1}(x)+X_{2}(y)+X_{3}(z),
\end{align*}
i.e.\ define $F_{i}(x)=C_{X_{i}(x)}$ where $C_{X}:H_{3}(\mk ) \to H_{3}(\mk
)$ is the commutator map
\begin{equation}
\label{eqn:com}
    C_{X}(H)=XH-HX.
\end{equation}
This satisfies equation~\eqref{eqn:Fone} and also
\begin{align}
\label{eqn:FonP}
    F_{i}(x)P_{i}(y) &= -2 \<x,y\> (e_{j}-e_{k}) \notag \\
    F_{i}(x)P_{j}(y) &= -P_k(\bar{y}\, \bar{x})  \\
    F_{i}(x)P_{k}(y) &= P_{j}(\bar{x}\, \bar{y}). \notag
\end{align}
It is a derivation of $H_{3}(\mk )$ by
virtue of the matrix identity \eqref{AHKJacobi}.

The Lie brackets of these derivations follow from another matrix
identity
\[
    [A,[B,H]]-[B,[A,H]]=[[A,B],H]+E(X,Y)H
\]
(see \eqref{ABHJacobi}). If
$A=X_{i}(x)$ and $B=X_{j}(y)$ we have $E(A,B)=0$ and
\begin{equation*}
    [X_{i}(x),X_{j}(y)]=X_{k}(\bar{y}\bar{x})
\end{equation*}
where $(i,j,k)$ is a cyclic permutation of $(1,2,3)$. This yields the
Lie bracket~\eqref{eqn:2}. If $X=X_{i}(x)$ and $Y=X_{i}(y)$, the matrix commutator
$Z=[X,Y]$ is diagonal with $z_{ii}=0$, $z_{jj}=y\bar{x}-x\bar{y}$ and
$z_{kk}=\bar{y}x-\bar{x}y$ ($i,j,k$ cyclic). Hence the action of the
commutator $[F_{i}(x),F_{i}(y)]=C_{Z}+E(X,Y)$ on $H_{3}(\mk )$ is
\begin{align*}
    [F_{i}(x),F_{i}(y)]e_{m}&= 0 \qqquad (m=i,j,k) \\
    [F_{i}(x),F_{i}(y)]P_{i}(w)&=
P_{i}\big( z_{jj}w-wz_{kk}-2[x,y,w]\big)\\
    &=4P_{i}(S_{xy}w) \quad \text{by eq. \eqref{S-def}.} \\
    [F_{i}(x),F_{i}(y)]P_{j}(w)&=P_{j}\big( z_{kk}w-2[x,y,w]\big) \\
    &=P_{j}(\bar{y}(xw)-\bar{x}(yw)) \\
    [F_{i}(x),F_{i}(y)]P_{k}(w) &= P_{k}\big( -wz_{jj}-2[x,y,w]\big) \\
    &= P_{k}\big((wx)\bar{y}-(wy)\bar{x}\big).
\end{align*}
Thus
\begin{align*}
    [F_{i}(x),F_{i}(y)]P_{i}(w) &=
P_{i}(Aw)=P_{i}(T_{i}w) \\
    [F_{i}(x),F_{i}(y)]P_{j}(w) &=
P_{j}(\overline{B}w)=P_{j}(T_j w) \\
    [F_{i}(x),F_{i}(y)]P_{k}(w) &=
P_{k}(\overline{C}w)=P_{k}(T_k w)
\end{align*}
where $T_{x,y}=(A,B,C)=(T_i,\overline{T_j},\overline{T_k})$, so that 
$(T_1,\overline{T_2},\overline{T_3})=\theta^{1-i}(T_{x,y})$. 
This establishes the Lie bracket~\eqref{eqn:3}.
\end{proof}
\end{theorem}

\begin{theorem}
\label{DerH3}
For any composition algebra $\mk$,
\begin{equation}
\label{eqn:derh3}
    \Der H_{3}(\mk ) = \Der \mk \ds A_{3}^{\prime}(\mk)
\end{equation}
in which $\Der \mk$ is a Lie subalgebra, the Lie brackets between $\Der
\mk $ and $A_{3}^{\prime}(\mk)$  are given by the elementwise action of
$\Der \mk$ on $3\times 3$ matrices over $\K$, and for $A,B\in A_3\pr(\K)$
\begin{equation*}
    [A,B]=(AB-BA)^{\prime}+\tfrac{1}{3}D(A,B)
\end{equation*}
where
\begin{equation*}
    D(A,B)=\sum_{ij} D_{a_{ij},b_{ji}} \quad \in \Der \mk ,
\end{equation*}
$a_{ij}$ and $b_{ij}$ being the matrix elements of $A$ and $B$.
\begin{proof}
By Lemma~\ref{TriK} and Theorem~\ref{theorem:banana}
\begin{equation}
\Der H_{3}(\mk ) = \Der \mk \ds 2\mk ^{\prime} \ds 3\mk .
\end{equation}
Identify $(a,b)+(x,y,z) \in 2\mk ^{\prime} \ds 3\mk$ with the traceless
antihermitian matrix
\begin{equation*}
    A= \begin{pmatrix} -a-b & -z & \bar{y} \\ \bar{z} & a & -x \\ -y
& \bar{x} & b \end{pmatrix} \quad \in A_{3}^{\prime}(\mk);
\end{equation*}
then the actions of $2\mk ^{\prime}$ and $3\mk $ on $H_{3}(\mk )$
defined in Theorem~\ref{theorem:banana} are together equivalent to the
commutator action $C_{A}$ defined by equation~\eqref{eqn:com}. By
the identity \eqref{CACB3},
\begin{equation*}
    [C_{A},C_{B}]= C_{(AB-BA)^{\prime}} + D(A,B),
\end{equation*}
so the bracket $[A,B]$ is as stated.
\end{proof}
\end{theorem}

It can be shown that this structure of $\Der\H_n(\K)$ persists for 
all $n$ if $\K$ is associative: 
\[
 \Der H_n(\K) \cong A\pr_n(\K)\dsum\Der\K.
\]
If $\K$ is not associative, however (i.e.\ $\K=\OO$ or $\widetilde{\OO})$,
the anticommutator algebra $H_n(\K)$ is not a Jordan algebra for $n>3$
and its derivation algebra collapses:
\[
 \Der H_n(\OO) = A\pr_n(\R)\oplus\Der\OO \cong \so(n)\oplus\g_2.
\]

\subsection{The Vinberg Construction}\label{Vinberg}
For this construction let $\mk _{1}$ and $\mk _{2}$ be composition
algebras, and let $\mk _{1}
\ox \mk _{2}$ be the tensor product algebra with multiplication
\begin{equation*}
(u_1\ox v_1)(u_2\ox v_2) = u_1v_1\ox u_2v_2
\end{equation*}
and conjugation
\begin{equation*}
\overline{u\ox v} = \overline{u}\ox\overline{v}.
\end{equation*}
Then the vector space
\begin{equation}
\label{eqn:vinberg}
    V_{3}(\mk _{1}, \mk _{2}) = A_{3}^{\prime}(\mk _{1}\ox \mk
_{2}) \ds \Der \mk _{1} \ds \Der \mk _{2}
\end{equation}
is clearly symmetric between $\mk_1$ and $\mk_2$. Vinberg showed that 
this is a Lie algebra when taken with the Lie
brackets defined by the statements:
\begin{enumerate}
\linespace\item $\Der \mk _{1} \oplus \Der \mk _{2}$ is a Lie subalgebra.
\linespace\item For $D \in \Der \mk _{1}\oplus \Der \mk _{2}$ and 
$A \in A_{3}^{\prime}(\mk _{1}\ox \mk_{2})$,
\begin{equation}\label{Vin1}
    [D,A]=D(A)
\end{equation}
where on the right-hand side $D$ acts elementwise on the matrix $A$.
\linespace\item For $A=(a_{ij})$, $B=(b_{ij}) \in A_{3}^{\prime}(\mk _{1}\ox \mk
_{2})$,
\begin{equation}\label{Vin2}
    [A,B] = (AB-BA)^{\prime}+\tfrac{1}{3}\sum _{ij} D_{a_{ij},b_{ji}}
\end{equation}
\end{enumerate}
where $D_{x,y}$ for $x,y\in\mk_1\ox\mk_2$ is defined by
\begin{equation*}
    D_{p\ox q, u\ox v} = \< p,u\> D_{q,v} +
\< q,v\> D_{p,u}.
\end{equation*}
We will now use the results of Section~\ref{Trider} to show that the Lie 
algebra $V_{3}(\mk _{1}, \mk _{2})$ is a Lie algebra isomorphic to the 
algebra $L_{3}(\mk _{1},\mk _{2})$ defined by the
Tits-Freudenthal construction, since no proof of this is readily
available.

\begin{theorem}\cite{Onishchik}
The Vinberg algebra defined above is a Lie algebra isomorphic to Tits's
magic square algebra $L_3(\K_1,\K_2)$.
\end{theorem}
\begin{proof}
The matrix part of the Vinberg vector space can be decomposed as
\[ 
A_3\pr(\K_1\ox\K_2) = A_3\pr(\R\ox\K_2) \dsum A_3\pr(\K_1\pr\ox\K_2).
\]
But $A_3\pr(\K_1\pr\ox\K_2)\cong\K_1\pr\ox H_3\pr(\K_2)$ (if $H$ is a
hermitian matrix over $\K_2$ and $a\in\K_1\pr$ is pure imaginary, then
$a\ox H$ is antihermitian over $\K_1\ox\K_2$), so
\[ 
A_3\pr(\K_1\ox\K_2)\cong A_3\pr(\K_2)\dsum\K_1\pr\ox H_3\pr(\K_2).
\]
Using Theorem~\ref{DerH3}, we can write the Tits vector space as
\[
L_3(\K_1,\K_2) = \Der\K_1\dsum\Der\K_2\dsum A_3\pr(\K_2)\dsum\K_1\pr\ox
H_3\pr(\K_2); 
\]
it is therefore isomorphic to $V_3(\K_1,\K_2)$. We will now prove that
this is a Lie algebra isomorphism by showing that
\begin{align*}
    [A_3\pr(\R\ox\K_2), A_3\pr(\R\ox\K_2)]_{\text{Vin}} 
&= [A_3\pr(\K_2),A_3\pr(\K_2)]_{\text{Tits}},\\
[A_3\pr(\R\ox\K_2), A_3\pr(\K_1\pr\ox H_3\pr(\K_2)]_{\text{Vin}}
&= [A_3\pr(\K_2), \K_1\pr\ox H_3\pr(\K_2)]_{\text{Tits}} \\
[A_3\pr(\K_1\pr \ox \K_2), A_3\pr(\K_1\pr\ox\K_2)]_{\text{Vin}}
&= [\K_1\pr\ox H_3\pr(\K_2),\K_1\pr\ox H_3\pr(\K_2)]_{\text{Tits}},
\end{align*}
where $[,]_{\text{Vin}}$ denotes the Lie brackets in the Vinberg
construction and $[,]_{\text{Tits}}$ denotes the Lie brackets in
the Tits construction.

\linespace
\textbf{1.} $[A_{3}^{\prime} (\R\ox\K_2),A_{3}\pr(\R\ox\K_2)]_{\text{Vin}}$.
In $V(\mk _{1},\mk _{2})$ the bracket is
\begin{equation*}
    [A,B]_{\text{Vin}} = (AB-BA)^{\prime} + \tfrac{1}{3}\sum _{ij}
D_{a_{ij},b_{ji}}
\end{equation*}
where $A, B \in A_{3}^{\prime} (\mk _{1}\ox \mr)$. In $L_3(\mk_1
,\mk_2)$ the matrices $A$ and $B$ are identified with elements of
$A_3^\prime(\mk_1)\subset\Der H_3(\mk_1)$, where their Lie bracket
$[A,B]_{\text{Tits}}$ is the same as the above by Theorem
\ref{DerH3}.

\linespace
\textbf{2.} $[A_{3}\pr(\R\ox\K_2), A_{3}\pr(\K_1\pr\ox\K_2)]_{\text{Vin}}$.
Let $A\in A_3\pr(\R\ox\K_2) = A_3\pr(\K_2)$ and $B\in
A_3\pr(\K_1\pr\ox\K_2)$; we may take $B = b\ox H$ with $b\in\K_1\pr$, $H\in
H_3\pr(\K_2)$. Then
\[ 
D_{a_{ij},b_{ji}} = \<1,b\>D_{a_{ij},h_{ji}} +
\<a_{ij},h_{ji}\>D_{1,b} =0
\]
and by Lemma \ref{trcomAH}, $\tr(AH-HA)=0$. Hence the Vinberg bracket is 
\[ 
[A,b\ox H]_{\text{Vin}} = b\ox(AH-HA) = [A,b\ox H]_{\text{Tits}}
\]
since the action of $A$ as an element of $\Der
H_3(\K_2)$ is $H\mapsto AH-HA$.

\linespace
\textbf{3.} $[A_3\pr(\K_1\pr\ox\K_2), A_3\pr(\K_1\pr\ox\K_2)]_{\text{Vin}}$. 
Let $A = a\ox H$, $B=b\ox H$ with $a,b\in\K_1\pr$ and $H,K\in H_3\pr
(\K_2)$. Then
\begin{align*}
[A,B]_{\text{Vin}} = &\, ab\ox HK - ba\ox KH - \third(ab\ox HK - ba\ox KH)\\
 &+ \third \sum_{ij}\(\<a,b\>D_{h_{ij},k_{ji}} + \<h_{ij},k_{ji}\>D_{a,b}\)\\
=& -\<a,b\>(HK-KH)\pr\\ 
= &\,\half[a,b]\ox(H*K)
- \<a,b\>D(H,K) + \third\<H,K\>D_{a,b}
\end{align*}
since $\Re(ab) = \Re(ba) = -\<a,b\>$ and $\Im(ab) = -\Im(ba) =
\half[a,b]$; also
\[ 
\sum D_{h_{ij},k_{ij}} = - \sum D_{h_{ij},\overline{k_{ij}}} = 
-\sum D_{h_{ij},k_{ji}} = - D(H,K)
\]
where $D(H,K)$ is defined in Theorem~\ref{DerH3}.

On the other hand,
\[ 
[a\ox H, b\ox K]_{\text{Tits}} = \third\<H,K\>D_{a,b} -
\<a,b\>[L_H,L_K] + \half[a,b]\ox(H*K).
\]
But the matrix identity \eqref{LHLK3} gives
\[ 
[L_H,L_K] = (HK-KH)\pr + \third D(H,K),
\]
from which it follows that  
\[
[a\ox H, b\ox K]_{\text{Tits}} = [a\ox H,b\ox K]_{\text{Vin}}.
\]
\end{proof}

\subsection{The Triality Construction}
A second clearly symmetric formulation of the magic square can be
given in terms of triality algebras.

\begin{theorem}
\label{theorem:peach}
For any two composition algebras $\mk _{1}$ and $\mk _{2}$,
\begin{equation}
\label{eqn:newL3}
L_{3}(\mk _{1}, \mk _{2}) = \Tri \mk _{1}\oplus \Tri \mk _{2} \ds 3\mk
_{1} \ox \mk _{2}
\end{equation}
in which $\Tri \mk _{1}\oplus \Tri \mk _{2}$ is a Lie subalgebra and the
other brackets are as follows. Define $F_i(x\ox y)\in
3\mk_1\ox\mk_2$ by
\[
F_{1}(x_{1}\ox x_{2}) + F_{2}(y_{1}\ox
y_{2}) + F_{3}(z_{1}\ox z_{2})=(x_{1}\ox x_{2},y_{1}\ox
y_{2},z_{1}\ox z_{2}).
\]
Then for $T_{\alpha }=(T_{\alpha 1},\overline{T}_{\alpha 2}, \overline{T}_{\alpha
3})\, \in \Tri \mk_\alpha$ and
$x_\alpha,y_\alpha,z_\alpha\in\mk_\alpha (\alpha = 1,2)$,
\begin{align}
    [T_{1},F_{i}(x_1\ox x_2)] &= F_{i}(T_{1i}x_1\ox x_2)
\label{eqn:th3.1} \\
    [T_{2},F_{i}(x_1\ox x_2)] &= F_{i}(x_1\ox T_{2i}x_2)
\label{eqn:th3.2} \\
[F_{i}(x_{1}\ox x_{2}) , F_{j}(y_{1}\ox y_{2})] &=
F_{k}(\overline{y}_{1}\overline{x}_{1}\ox \overline{y}_{2}\overline{x}_{2})
\label{eqn:th3.3}
\end{align}
where $(i,j,k)$ is a cyclic permutation of $(1,2,3)$; and
\begin{multline}
\label{eqn:th3.4}
[F_{i}(x_{1}\ox x_{2}) , F_{i}(y_{1}\ox y_{2})] = \<
x_{2},y_{2}\> \theta ^{1-i}T_{x_{1}y_{1}}+\<
x_{1},y_{1}\> \theta ^{1-i}T_{x_{2}y_{2}} \\
\hspace{4.5cm} \in \Tri \mk_{1} \oplus \Tri \mk_{2}
\end{multline}
where $\theta$ is the automorphism of Lemma \ref{lemma:splodge}.
\begin{proof}
The vector space structure~\eqref{Titsspace} of $L_{3}(\mk _{1},\mk_{2})$
can be written, using Theorem~\ref{theorem:banana} and Lemma~\ref{TriK}, as
\begin{align*}
    L_{3}(\mk _{1},\mk_{2}) &= \Der H_{3}(\mk _{1}) \ds H_{3}^{\prime
}(\mk _{1})\ox \mk _{2}^{\prime} \ds \Der \mk_{2} \\
    &= (\Tri \mk_1 \ds 3\mk _{1})\ds (2\mk_{2}^{\prime} \ds
3\mk_{1}\ox \mk _{2}^{\prime}) \ds \Der \mk_{2} \\
    &= \Tri \mk _{1} \ds (\Der \mk _{2} \ds 2\mk _{2}^{\prime}) \ds
(3\mk_{1}\ox \mk_{2}^{\prime} \ds 3\mk _{1}) \\
    &\cong \Tri \mk _{1} \ds \Tri \mk _{2} \ds 3\mk _{1}\ox \mk
_{2}.
\end{align*}

We need to consider the following five subspaces of
$L_{3}(\mk _{1},\mk _{2})$:
\begin{enumerate}
\linespace\item $\Tri \mk \subset \Der H_{3}(\mk _{1})$ contains elements
$T=(T_{1},\overline{T}_{2},\overline{T}_{3})$ acting on
$H_{3}^{\prime}(\mk _{1})$ as in Theorem~\ref{theorem:banana}:
\begin{equation*}
Te_{i}=0, \quad TP_{i}(x)=P_{i}(T_{i}x) \quad (x\in \mk ; i = 1,2,3).
\end{equation*}
\item $3\mk _{1}$ is the subspace of $\Der H_{3}(\mk _{1})$ containing
the elements $F_{i}(x)$ defined in Theorem~\ref{theorem:banana}; these
will be identified with the elements $F_{i}(x\ox 1)\in 3\mk
_{1}\ox \mk _{2}$.
\linespace\item $2\mk _{2}^{\prime}$ is the subspace $\Delta \ox \mk
_{2}^{\prime}$ of $H_{3}(\mk _{1})\ox \mk _{2}^{\prime}$, where
$\Delta \subset H_{3}^{\prime}(\mk _{1})$ is the subspace of
real, diagonal, traceless matrices and is identified with a subspace of
$\Tri \mk_2 $ as described in Lemma~\ref{TriK}. We will regard $2\mk
_{2}^{\prime}$ as a subspace of $3\mk_{2}^{\prime}$, namely
\begin{equation*}
    2\mk _{2}^{\prime}= \{ (a_{1},a_{2},a_{3})\in 3\mk _{2}^{\prime}
: a_{1}+a_{2}+a_{3}=0 \}
\end{equation*}
and identify $\mathbf{a}=(a_{1},a_{2},a_{3})$ with the $3\times 3$
matrix
\begin{equation*}
    \Delta (\mathbf{a})= \begin{pmatrix} a_{1} & 0 & 0 \\ 0 & a_{2}
& 0 \\ 0 & 0 & a_{3} \end{pmatrix} \in H_{3}^{\prime}(\mk _{1})\ox
\mk _{2}^{\prime}
\end{equation*}
in the Tits description, and on the other hand with the triality
$T(\mathbf{a})=(T_{1},\overline{T}_{2},\overline{T}_{3})$ where
$T_{i}=L_{a_{j}}-R_{a_{k}}$.

\linespace\item $3\mk _{1}\ox \mk _{2}^{\prime}$ is the subspace of
$H_{3}(\mk _{1}\ox \mk _{2}^{\prime})$ spanned by elements
$P_{i}(x)\ox a$ $(i=1,2,3 : x \in \mk _{1}, a \in \mk
_{2}^{\prime})$; in the triality description it is a subspace of 
$3\mk _{1}\ox \mk _{2}$ in the obvious way.

\linespace\item $\Der \mk _{2}$ is a subspace of $\Tri \mk _{2}$, a derivation $D$
being identified with $(D,D,D) \in \Tri \mk _{2}$.
\end{enumerate}

The proof is completed by verifying that the Lie brackets defined by
Tits (eq.~\ref{Titsbrac}) coincide with those in the statement of the
theorem. The above decomposition of $L_3(\K_1,\K_2)$ gives fifteen types
of bracket to examine; for each of them the verification is
straightforward \cite{Barton}.
\end{proof}
\end{theorem}

The isomorphism between the Vinberg construction and the triality
construction is easy to see directly at the vector space level: using
Lemma~\ref{TriK}, 
\begin{align*}
    L_{3}(\mk_{1},\mk_{2}) &= \Tri\K_1\dsum\Tri\K_2\dsum 3\K_1\ox\K_2\\
&=\Der\K_1\ds 2\K_2\pr\ds\Der\K_2\ds2\K_2\pr\dsum 3\K_1\ox\K_2\\
&=\Der\K_1\dsum 2\K_1\pr\ox\R\dsum\Der\K_2\dsum 2\R\ox\K_2\pr\ds
3\K_1\ox\K_2\\
&= A_{3}(\K_1\ox \K_2) \ds \Der \K_1\ds \Der \K_2.
\end{align*}

Thus both ways of understanding the symmetry of the $3\times 3$ magic
square reduce to being different ways of looking at the same underlying
vector space, which is an extension of the vector space of
antisymmetric $3\times 3$ matrices over $\mk_1\ox\mk_2$. The Lie
algebras of the square can therefore be understood as analogues of
$\su(3)$ with the complex numbers replaced by $\mk_1\ox\mk_2$.

\section{The rows of the magic square}\label{rows}

In this section we will examine the non-symmetric magic square obtained
by taking $\K_1$ to range over the split composition algebras $\R$,
$\sC$, $\sH$ and $\sO$. According to Theorem~\ref{Jmatrix}, the first three
rows contain the derivation, structure and conformal algebras of the
Jordan algebras $H_3(\K)$, which we have defined to be the
generalisations of the Lie algebras of antihermitian traceless $3\times
3$ matrices, all traceless $3\times 3$ matrices, and symplectic
$6\times 6$ matrices:
\begin{align*}
L_3(\R,\K) &\cong \Der H_3(\K) = \sa (3,\K),\\
L_3(\sC,\K) &\cong \Str\pr H_3(\K) = \sl (3,\K),\\
L_3(\sH,\K) &\cong \Con H_3(\K) = \sp(6,\K).
\end{align*}
We will now determine the precise composition of these algebras in terms
of matrices over $\K$.

\begin{theorem}\label{3x3mat}

\emph{\textbf{(a)}}\upline
\begin{equation}
\sa(3,\K) = A_3\pr(\K) \dsum \Der\K; \label{samat}
\end{equation}

\emph{\textbf{(b)}}\upline
\begin{equation}
\sl(3,\K) = L_3\pr(\K) \dsum \Der\K; \label{slmat}
\end{equation}

\emph{\textbf{(c)}}\upline
\begin{equation}
\sp(6,\K) = Q_6\pr(\K) \dsum \Der\K. \label{spmat}
\end{equation}

In each case the Lie brackets are defined as follows:
\begin{enumerate}
\item $\Der\K$ is a Lie subalgebra;
\item The brackets between $\Der\K$ and the other summand are given by
the elementwise action of $\Der\K$ on matrices over $\K$;
\item The brackets between two matrices in the first summand are
\be \label{mat3brac}
[X,Y] = (XY-YX)\pr + \frac{1}{n}D(X,Y) 
\end{equation}
where $n$ (= 3 or 6) is the size of the matrix and $D(X,Y)$ is defined
in \eqref{Dmatdef}.
\end{enumerate}

\begin{proof}
\textbf{(a)} This is Theorem~\ref{DerH3}.

\linespace
\textbf{(b)} The vector space of $\sl(3,\K)$ is
\begin{align*}
\sl(3,\K) &= \Str\pr H_3(\K) = \Der H_3(\K) \dsum H_3\pr(\K)\\
&= \Der\K\dsum A_3\pr(\K)\dsum H_3\pr(\K)\\
&= \Der\K\dsum M_3\pr(\K)
\end{align*}
For $A,B\in A_3\pr(\K)$ the Lie bracket is that of $\Der H_3(\K)$, which
is \eqref{mat3brac}. For $A\in A_3\pr(\K)$,
$H\in H_3\pr(\K)$ the bracket is given by the action of $X$ as an
element of $\Der H_3(\K)$ on $H$, which according to Theorem \ref{DerH3} is
\be \label{XHcom}
[A,H] = AH - HA.
\end{equation}
Now by Lemma \ref{trcomAH}, $\tr(AH-HA)=0$ and
$D(A,H)=0$; hence \eqref{XHcom} is the same as \eqref{mat3brac}.  

Finally, for $H,K\in H_3\pr(\K)$ the $\Str\pr H_3(\K)$ bracket is 
\[
[H,K] = L_H L_K - L_K L_H \in \Der H_3(\K)
\]
and in the Appendix it is shown that this commutator appears in the 
decomposition $\Der H_3(\K)$ as
\[ [L_H,L_K] = (HK-KH)\pr + \third D(H,K).
\]

Note that the action of Str$H_3(\K)$ on $H_3(\K)$ is as follows. The
subalgebra $\Der\K$ acts elementwise, while according to
\eqref{DerH3} the matrix part of $\Der H_3(\K)$ acts by
\[
H \mapsto AH - HA \qqquad (A\in A_3\pr(\K)).
\]
The remaining matrix subspace $H_3\pr(\K)$ acts by translations in the
Jordan algebra $H_3(\K)$:
\[
H \mapsto KH + HK \qqquad (K \in H_3\pr(\K)).
\]
Hence the action of the matrix part of Str$\pr H_3(\K)$ is
\be \label{StrH3act}
H \mapsto XH + HX^\dagger \qqquad (X\in M_3\pr (\K)).
\end{equation}

\linespace
\noindent\textbf{(c)} The vector space of $\sp(6,\K)$ is
\begin{align}
\sp(6,\K) &= \Con H_3(\K) = \Str H_3(\K) \dsum 2H_3(\K)\notag\\
&= \so(\K\pr) \dsum M_3\pr(\K) \dsum \R \dsum 2H_3(\K).\label{spvec}
\end{align}
On the other hand, a $6\times 6$ matrix $X$ belongs to $Q_6\pr(\K)$ if
and only if
\begin{align*}
X^\dagger J + JX &= 0 \quad \text{and} \quad \tr X = 0\\
\iff X &= \begin{pmatrix} A & B\\C & -A^\dagger \end{pmatrix} 
\quad \text{ with } B,C\in H_3(\K)\\
 \text{ and }\quad &A \in M_3(\K), \quad\Im(\tr A) = 0, \quad \text{ so } 
\quad A\in M_3\pr(\K)\dsum\R.
\end{align*}
Thus
\[ 
Q_6\pr(\K) \cong M_3\pr(\K) \dsum \R \dsum 2H_3(\K),
\]
the summand $\R$ representing $\Re(\tr A)$, so the vector space
structure of $\sp(6,\K)$ is as stated in \text{(c)}.

To examine the Lie brackets, we write \eqref{spmat} as 
\[
\sp(4,\K) \cong \sl(3,\K) \dsum \R \dsum 2H_3(\K).
\]
An element $A$ of the matrix part of $\sl(3,\K)$ corresponds in
$\sp(6,\K)$ to the matrix
$\widehat{A}=\begin{pmatrix}A&0\\0&-A^\dagger\end{pmatrix}$. Since
Str$H_3(\K)$ is a Lie subalgebra of Con$H_3(\K)$, the Lie bracket of two
such elements in $\sp(6,\K)$ is given by 
\[
[\widehat{A},\widehat{B}] = (AB-BA-\third tI_3)\widehat{\phantom{a}} + 
\third D(A,B)
\]
where $t = \tr(AB-BA)$, which is purely imaginary, being a sum of
commutators in $\K$. Hence 
\begin{align*}
(AB-BA-\third tI_2)\widehat{\phantom{a}} &= \widehat{A}\widehat{B} - 
\widehat{B}\widehat{A} - \third tI_6\\
&= \widehat{A}\widehat{B} - \widehat{B}\widehat{A} - 
\tfrac{1}{6}\tr(\widehat{A}\widehat{B} - \widehat{B}\widehat{A})I_6\\
&= (\widehat{A}\widehat{B} - \widehat{B}\widehat{A})\pr.
\end{align*}
Also $D(A,B)=2D(\widehat{A},\widehat{B})$, so \eqref{mat3brac} holds in
$\sp(6,\K)$ for elements $X,Y$ of the form $\widehat{A}$.

For $X\in M_3\pr(\K)$ and $Y\in\R$ or $X,Y\in\R$, both sides of
\eqref{mat3brac} are zero.

For $X\in M_3\pr(\K)\dsum\R$ and $Y\in 2H_3(\K)$, i.e.
\[ 
X=\begin{pmatrix}A&0\\0&-A^\dagger\end{pmatrix}, \qquad
Y=\begin{pmatrix}0&B\\C&0\end{pmatrix}
\]
with $A\in M_3(\K)$, $B,C\in H_3(\K)$, the Lie bracket in $\sp(6,\K)$
is given by the direct sum of the action of $\Str H_3(\K)$ on $H_3(\K)$
and its transform by the involution * of Section~\ref{notation}. The
action is given by \eqref{StrH3act}, so $B\mapsto AB + BA^\dagger$,
while the effect of the involution is to change the sign of the
hermitian part of $A$, so $C\mapsto -A^\dagger C - CA$. Thus
\[
[X,Y] = \begin{pmatrix}0 & AB+BA^\dagger\\
-A^\dagger C - CA & 0\end{pmatrix} = XY - YX.
\]
Clearly $\tr(XY-YX) =0$ and $D(X,Y)=0$, so this is the same as
\eqref{mat3brac}.

Finally, for $X,Y\in 2H_3(\K)$, say
\[
X=\begin{pmatrix}0 & H\\K & 0\end{pmatrix} \quad \text{ and } \quad
Y=\begin{pmatrix}0 & B\\C & 0\end{pmatrix},
\]
the bracket is given by (\ref{Conbrac2} - \ref{Conbrac3}), i.e.
\[
[X,Y] = \half\(L_{H\cdot C} + L_{K\cdot B} + [L_H,L_C] + [L_K, L_B]\).
\]
The first two terms on the right-hand side form an element of
$\Str H_3(\K)$ which corresponds in $\sp(6,\K)$ to the matrix
\[
\frac{1}{2}\begin{pmatrix} HC+CH+KB+BK & 0 \\ 0 & -HC-CH-KB-BK
\end{pmatrix}
\]
while the second pair of terms forms an element of $\Der H_3(\K)$
corresponding, according to \eqref{mat3brac}, to the sum of the matrix
\[
\frac{1}{2}\begin{pmatrix} (HC-CH)\pr + (KB-BK)\pr & 0 \\
0 & (HC-CH)\pr + (KB-BK)\pr \end{pmatrix}
\]
and the $\Der\K$ element
\[
\sixth D(H,C) + \sixth D(H,C) + \sixth D(K,B) = \sixth D(X,Y).
\]
Hence
\begin{align*}
[X,Y] &= \begin{pmatrix}HC+KB & 0 \\ 0 & -CH-BK \end{pmatrix}\\
&\qquad -\sixth\tr(HC-CH+KB+BK)I_6 + \sixth D(X,Y)\\
&= (XY-YX)\pr + \sixth D(X,Y)
\end{align*}
as asserted in \text{(c)}.
\end{proof}
\end{theorem}

This description of the rows of the magic square was given a geometrical
interpretation by Freudenthal \cite{Freudenthal65}. A Lie algebra of
$3\times 3$ matrices corresponds to a Lie group of linear
transformations of a 3-dimensional vector space, or projective
transformations of a plane. The Lie group corresponding to
$\sa(3,\K)$ preserves a hermitian form in the vector space or a
\emph{polarity} in the projective plane, i.e.\ a correspondence between
points and lines. This defines the four (real, complex, quaternionic and
octonionic) \emph{elliptic} geometries, in which there is just one class
of primitive geometric objects, the points, with a relation of polarity
between them (inherited from orthogonality of lines in the vector space).
The special linear Lie algebras $\sl(3,\K)$ correspond to the
transformation groups of the four projective geometries, in which there
are two primitive geometric objects, points and lines, with no relations
between points and points or lines and lines but a relation of incidence
between points and lines. The third row of the magic square, containing Lie
algebras $\sp(6,\K)$, yields the transformation groups of
five-dimensional symplectic geometries, whose primitive geometric
objects are points, lines and planes. Freudenthal completed this
geometrical schema to incorporate the last row of the magic square by
defining \emph{metasymplectic} geometries, which have a fourth type of
primitive object, the \emph{symplecta}. In metasymplectic geometry
points can be \emph{joined} (contained in a line, which is unique if the
points are distinct), \emph{interwoven} (contained in a plane, unique if
the points are not joined) or \emph{hinged} (contained in a symplecton,
unique if the points are not interwoven).

\section{Magic Squares of $n\times n$ Matrices}\label{nxn}

According to Theorem \ref{thm:Tits}, Tits's construction 
(\ref{Titsspace}--\ref{Titsbrac}) yields a Lie algebra for any Jordan
algebra if the composition algebra $\K_2$ is associative. Hence for
$\K_2 = \R,\C,\H$ and their split versions we obtain a Lie algebra
$L_n(\K_1,\K_2)$ for any $n>3$ by taking $\J=H_n(\K_2)$ in $L(\K_1,\J)$
(the case $n=2$, which will be examined in section~\ref{2mat}, lends
itself naturally to a slightly different construction). The proof of
Vinberg's model is valid for any size of matrix, so we have

\begin{theorem}\label{thm:ntimesn} Let $\K_1$ and $\K_2$ be associative
composition algebras over $\R$, and let $L_n(\K_1,\K_2)$ be the Lie
algebra obtained by Tits's construction
\emph{(\ref{Titsspace}--\ref{Titsbrac})} with 
$\K=\K_1$ and $\J=H_n(\K_2)$. Then
\[
 L_n(\K_1,\K_2) = A_n\pr(\K_1\otimes\K_2) \dsum \Der\K_1 \dsum \Der\K_2
\]
with brackets as in $V_3(\K_1,\K_2)$ (section \ref{Vinberg}).
\end{theorem}

By Theorem~\ref{Jmatrix}, the Lie algebras $\sa(n,\K)$,
$\sl(n,\K)$ and $\sp(2n,\K)$ can now be identified for associative $\K$
as 
\begin{align}
\sa(n,\K) &= \Der H_n(\K) = L_n(\R,\K) = A_n\pr(\K) \dsum \Der\K,\\
\sl(n,\K) &= \Str\pr H_n(\K) = L_n(\mcs,\K) = M_n\pr(\K) \dsum
\Der\K,\\
\sp(2n,\K) &= \Con H_n(\K) = L_n(\mhs,\K) = Q_n\pr(\K) \dsum \Der\K.
\end{align}

\subsection{The Santander-Herranz Construction}

Vinberg's approach to the magic square
is extended to general dimensions $n$ by Santander and
Herranz in their construction of
`Cayley-Klein' (CK) algebras. This starts from
a $(2N+2)\times(2N+1)$ matrix $I_{\omega}=\diag (1,
\omega _{01}, \omega _{02}, \dots , \omega _{0N})$ depending on $N+1$ fixed
non-zero parameters $\omega _{i}$, with $\omega
_{0a}=\omega _{0}\omega_{1}\dots \omega _{a}$. Let $\mathbb{I}_{\omega} =
\left( \begin{smallmatrix} 0 & I_{\omega} \\ -I_{\omega} & 0
\end{smallmatrix} \right)$. A
matrix $X$ is defined to be $G$-antihermitian if $X^{\dagger}G+GX=0$.
Santander and Herranz define three series of classical
CK-algebras:
\begin{enumerate}
\item The special antihermitian CK-algebra, $\sa _{\omega _{1} \dots
\omega _{N}}(N+1, \mk)$. This is the Lie algebra of
$I_{\omega}$-antihermitian matrices, $X$, over $\mk$ if $\mk =
\mr$ or $\mh$, or the subalgebra of traceless matrices if $\mk =\mc$.
\linespace\item The special linear CK-algebra, $\sl _{\omega _{1} \dots
\omega _{N}}(N+1,\mk )$. This is the Lie algebra of all matrices $X\in
\K^{(N+1)\times(N+1)}$ with $\tr X = 0$ if $\mk = \mr$ or $\mc$ and $\Re (\tr X)=0$
if $\mk=\mh$.
\linespace\item The special symplectic CK-algebra, $\mathfrak{sn} _{\omega _{1} \dots
\omega _{N}}(2(N+1),\mk)$. This is the Lie algebra of all
$\mathbb{I}_{\omega}$-antihermitian matrices over $\mk$ if $\mk = \mr$
or $\mh$ and the subalgebra of matrices with zero trace if $\mk =\mc$.
\end{enumerate}
For $N=1,2$ these definitions can be extended to include $\mk = \mo$ by
adding the derivations of $\mo$ in each case.
A fourth CK-algebra can also be added, the metasymplectic CK-algebra,
$\mathfrak{mn} (N+1,\mk)$ based on the definition of the metasymplectic geometry
given in~\cite{Freudenthal65}.

Now define the set of matrices
\begin{equation*}
J_{ab}=\begin{pmatrix}  & \vdots & &\vdots & \\
\cdot & \cdot & \cdots & -\omega_{ab} & \cdots \\
 & \vdots & & \vdots & \\
\cdot & 1 & \cdots & \cdot & \cdots \\
 & \vdots & & \vdots
\end{pmatrix},  \quad M_{ab}=\begin{pmatrix}  & \vdots & &\vdots & \\
\cdot & \cdot & \cdots & \omega_{ab} & \cdots \\
 & \vdots & & \vdots & \\
\cdot & 1 & \cdots & \cdot & \cdots \\
 & \vdots & & \vdots
\end{pmatrix}
\end{equation*}
and
\begin{equation*}
H_{m} = \begin{pmatrix} 1 &  & \cdot & \\
 & & \vdots & \\
\cdot & \cdots & 1 & \cdots \\
 & & \vdots & \end{pmatrix}, \quad E_{0}=\begin{pmatrix} 1 & \cdots \\
\vdots & 0 & \end{pmatrix},
\end{equation*}
where $a,b=0,1,\dots,N$ with the condition that $a<b$; $m=1,\dots,N$; and
matrix indices run over the range $0,\dots ,N$. Further if $X$ is one of
these matrices then define $X^{i}=e_{i}X$ and
\begin{equation*}
\mathbb{X} = \begin{pmatrix} X & 0 \\ 0 & X \end{pmatrix},\;
\mathbb{X}_{1} =\begin{pmatrix} 0 & X \\ -X & 0 \end{pmatrix},\;
\mathbb{X}_{2} =\begin{pmatrix} 0 & X \\ X & 0 \end{pmatrix},\;
\mathbb{X}_{3} =\begin{pmatrix} X & 0 \\ 0 & -X \end{pmatrix}.
\end{equation*}
Note that there is an isomorphism $J\mapsto
\mathbb{J}_{ab}$, $M_{ab}\mapsto \mathbb{M}_{ab;2}$, $M_{ab}^{1} \mapsto
\mathbb{M}_{ab}^{1}$, $M_{ab}^{2} \mapsto
\mathbb{M}_{ab}^{2}$.

The first three rows and columns of the Tits-Freudenthal magic square
can now be generalised to the $(N+1)$-dimensional case using the three
CK-algebra series as follows
\begin{center}
\vspace{0.5cm}
\begin{tabular}{|c|c|c|c|}
\hline
\multicolumn{1}{|c|}{Lie Algebra} &\multicolumn{3}{|c|}{Lie span of the
generators} \\ \cline{2-4}
 & $\mr$ & $\mc$ & $\mh$ \\
\hline
$\sa _{\omega _{1}\dots \omega _{N}}(N+1,\mk)$ & $J_{ab}$ &
$J_{ab}$,$M_{a,b}^{1}$ & $J_{ab}$,$M_{a,b}^{1}$,$M_{ab}^{2}$ \\
\hline
$\sl _{\omega _{1}\dots \omega _{N}}(N+1,\mk)$ & $J_{ab}$,$M_{ab}$ &
$J_{ab}$,$M_{ab}$,$M_{a,b}^{1}$ & $J_{ab}$,$M_{ab}$,$M_{a,b}^{1}$,$M_{ab}^{2}$
\\
\hline
$\mathfrak{sn} _{\omega _{1}\dots \omega _{N}}(N+1,\mk)$ &  $\mj
_{ab}$,$\mathbb{M}_{ab;1}$,$\mathbb{M}_{ab;2}$ & $\mj
_{ab}$,$\mathbb{M}_{ab;1}$,$\mathbb{M}_{ab;2}$,$\mathbb{M}_{ab}^{1}$ & $\mj
_{ab}$,$\mathbb{M}_{ab;1}$,$\mathbb{M}_{ab;2}$,$\mathbb{M}_{ab}^{1}$,
$\mathbb{M}_{ab}^{2}$ \\
\hline
\end{tabular}
\vspace{0.5cm}
\end{center}
Then the symmetry of the $(N+1)$ dimensional magic square (and
consequently of the $3\times 3$ magic square) can be
explained as follows.

Each algebra is a subalgebra of all the algebras to its right and below
it: as we move from left to
right and from top to bottom across the square, in each step the same new
generators appear. Explicitly, moving from the top algebra ($\sa$) to
the bottom ($\mathfrak{sn}$), in each column $M_{ab}$ appears in the
first step ($\sa \rightarrow \sl$) and $\mathbb{M}_{ab;1}$ appears in
the second ($\sl \rightarrow \mathfrak{sn}$). Similarly, moving from left
to right, $M_{ab}^{1}$ is the additional generator after the first step
and $M_{ab}^{2}$ is the additional generator after the second.

In more recent work Santander~\cite{Santander2} has gone on to define the tensor algebra
$\sa _{\omega _{1}\dots \omega _{N}}(N+1,\mk _{1}\ox \mk _{2})$,
an extension of the Vinberg construction which includes all simple
Lie algebras, i.e.\ \textit{any} simple Lie algebra can be written in the
form $\sa _{\omega _{1}\dots \omega _{N}}(N+1,\mk _{1}\ox \mk _{2})$
for an appropriate choice of $\omega _{i}$, $N$, $\mk _{1}$ and $\mk _{2}$.
Explicitly this is the algebra of $(N+1)\times (N+1)$ matrices with
entries in $\mk _{1} \ox \mk_{2}$ and the derivations of $\mk _{1}$
and $\mk _{2}$.

Thus we have a second way of approaching an explanation of the symmetry
of the magic square and indeed a classification of all simple Lie algebras in
terms of matrices with entries in the division algebras.

\section{Maximal Compact Subalgebras}\label{maxcomp}

We now turn to the question of identifying the Lie algebras of
Theorem~\ref{3x3mat} in the standard list of real forms of complex
semisimple Lie algebras. The split magic square $L_3(\sK_1,\K_2)$
contains real forms of the complex Lie algebras $L_3(\K_1,\K_2)$ which
are identified in Table 1; we will establish this identification by finding
the maximal compact subalgebras.

Recall that a semi-simple Lie algebra over $\R$ is called \emph{compact} if it has a
negative-definite Killing form. A non-compact real form $\g$ of a semi-simple
complex Lie algebra $L$ has a
maximal compact subalgebra $\n$ with an orthogonal complementary
subspace $\p$ such that $\g = \n \ds \p $ and the brackets
\begin{align}
    [\n ,\n ] & \subseteq \n \notag \\
    [\n ,\p ] & \subseteq \p \label{comp-decomp}\\
    [\p ,\p ] & \subseteq \n \notag
\end{align}
(see, for example, ~\cite{Gilmore}), from which it follows that 
$\<\n ,\p \> = 0$ where $\<,\>$ is the Killing
form of $L$. There exists an involutive automorphism $\sigma : \g
\rightarrow \g$ such that $\n $ and $\p $ are eigenspaces of $\sigma $ with
eigenvalues $+1$ and $-1$ respectively. A compact real form, $\g
^{\prime}$, of $L$ will also contain $\n$ as a compact subalgebra of $\g
^{\prime}$ but clearly in this case the maximal compact subalgebra will be $\g
^{\prime}$ itself. We can obtain $\g ^{\prime}$ from $\g$ by keeping the
same brackets in $[\n ,\n ]$ and $[\n ,\p ]$ but multiplying the
brackets in $[\p ,\p ]$ by $-1$, i.e.\ by performing the \emph{Weyl
unitary trick} (putting $\g
^{\prime}= \n \ds i\p $).

We will use the following method to identify the maximal compact
subalgebras in $L_3(\sK_1,\K_2)$. It is known that
$L_{3}(\mk _{1},\mk _{2})$ gives a compact real form of each Lie algebra
 (from, for example~\cite{JacobsonELA}).
Thus if  $L_{3}(\sK_{1},\K_{2})$ shares a common subalgebra with
$L_{3}(\sK_{1},\K_{2})$, say $\n $, where
\begin{align*}
    L_{3}(\mk _{1},\mk _{2}) &= \n \ds \p _{1} \\
    L_{3}(\mk _{1},\mks _{2})&= \n \ds \p _{2},
\end{align*}
and the brackets in $[\n , \p _{1}]$ are the same as those in $[\n ,\p
_{2}]$ but the brackets in $[\p _{1}, \p _{1}]$ are $-1$ times the
equivalent brackets in $[\p _{2},\p _{2}]$, then $\n $ will be the
maximal compact subalgebra of $L_{3}(\mk _{1},\mks _{2})$ and $\p _{2}$
will be its orthogonal complementary subspace. We will see that
this sign change in the
brackets reflects precisely the change in sign in the Cayley-Dickson
process~\cite{Schafer66} when moving from the division algebra to the corresponding
split composition algebra.

First we consider the relation between $\Der\K$ and $\Der\sK$ where $\K$
is a division algebra. Both
$\K$ and $\sK$ can be obtained by the Cayley-Dickson
process~\cite{Schafer66} from a positive-definite composition algebra $\F$; they are both of
the form $F_\epsilon^{\,2} = \F \dsum l\F$, where $l$ is the new imaginary
unit and the multiplication is given by
\begin{align}
x(ly) &= l(\overline{x}y)\notag\\
(lx)y &= l(yx)\label{C-D}\\
(lx)(ly) &= \epsilon y\overline{x}\notag
\end{align}
where $\epsilon = -1$ for $\K$ and $\epsilon = 1$ for $\sK$.

A derivation $D$ of $\F_\epsilon^{\,2}$ can be specified by giving its
action on $\F^{\,\prime}$ (since $D(1)=0$) and by specifying $D(l)$. Thus
each derivation $D$ of $\F$ can be extended to a derivation
$\overline{D}$ of $\F_\epsilon^{\,2}$ by defining $Dl=0$. We also define
derivations $E_a, F_a$ for each $a\in \F^{\,\prime}$, and $G_S$ for each
symmetric linear map $S:\F^{\,\prime} \to \F^{\,\prime}$, as follows:
For $x\in\F$ and $a,b\in\F^{\,\prime}$,
\begin{align}\label{deriv}
E_a x &= 0 \quad (x\in\F), &E_a l &= la;\notag\\
F_a b = l(ab &-\<a,b\>) \quad (b\in\F^{\,\prime}), &F_al &= -2\epsilon a\\
&&\text{so that }F_a(lb) &= -\half\epsilon[a,b];\notag\\
G_Sa &= l(Sa)\quad (a\in\F^{\,\prime}), &G_Sl &= 0.\notag
\end{align}

\begin{theorem}\label{thm:deriv}
$\Der(\F_\epsilon^{\,2})$ is spanned by
$\overline{D}$ $(D\in\Der\F),\;E_a,\;F_a\,(a\in\F^{\,\prime})$ and $G_S$
where $S:\F^{\,\prime}\to\F^{\,\prime}$ is symmetric and traceless if $\F=\H$.
The Lie brackets are given by
\[ [D,E_a] = E_{Da},\qquad [D,F_a] = F_{Da}, \]
\[ [D,G_S] = G_{[D,S]}, \]
\[ [E_a, E_b] = -E_{[a,b]} \]
\[ [E_a, F_b] = \quarter F_{[a,b]} - \tfrac{3}{2}G_{S(a,b)}
+\(1-\tfrac{3}{m}\)\<a,b\>G_{\id} \]
\emph{(}where $m=\dim \F^\prime$, $S(a,b)$ is the traceless symmetric map
\[
S(a,b)c = \<a,c\>b + \<b,c\>a - \frac{2}{m}\<a,b\>c, 
\]
and $G_{\id}$ is given by (\ref{deriv}) when $S$ is the identity map on
$\F^\prime$\emph{);}
\[ 
[E_a, G_S] = \half F_{Sa} - \quarter G_{[D_a,S]} 
\]
where $D_a \in \Der\F$ is the inner derivation $D_a(x) = [a,x]$;
\[ [F_a, F_b] = -\quarter\epsilon D_{[a,b]} -2\epsilon E_{[a,b]},\]
\[ [F_a, G_S] = \half\epsilon D_{Sa} + 2\epsilon  E_{Sa},\]
\[ [G_S,G_T] = \epsilon [S,T].\]
\end{theorem}

\begin{proof} First we note that any derivation of an algebra must
annihilate the identity of the algebra. Let $D$ be a derivation of
$\F_\epsilon^{\,2}$ satisfying $Dl=0$. Then $D$ is determined by its action
on $a\in\F^{\,\prime}$. Write 
\[
 Da = Ta +l(Sa)
\]
where $T$ and $S$ are maps from $\F^{\,\prime}$ to $\F$. Then the derivation
condition applied to the relations \eqref{C-D} requires $T$ to be a
derivation of $\F$ and $S$ to be a map from $\F^{\,\prime}$ to $\F^{\,\prime}$
satisfying 
\begin{align*} 
S[a,b] &= -2a(Sb) + 2b(Sa)\\
&= 2(Sb)a - 2(Sa)b.
\end{align*}
Since $ab+ba = -2\<a,b\>$ where $\<,\>$ is the inner product on $\F$,
this yields 
\[
 \<Sa,b\> = \<a,Sb\>,
\]
i.e.\ $S$ is a symmetric operator on $\F^{\,\prime}$, and 
\[
 S[a,b] = -[Sa,b] - [a,Sb].
\]
If $\F=\H$ there is an identity
\[
 S[a,b] + [Sa,b] + [a,Sb] = (\tr S)[a,b]
\]
(a version of $\epsilon_{mjk}\delta_{in} + \epsilon_{imk}\delta_{jn} +
\epsilon_{ijm}\delta_{kn} = \epsilon_{ijk}\delta_{mn}$),
so that in this case $\tr S = 0$.

Thus 
\[
 Dl = 0 \implies D = T + G_S \;\text{ with }\;T\in\Der\F
\]
where $S:\F^{\,\prime}\rightarrow\F^{\,\prime}$ is symmetric and traceless if $\F=\H$. To
show that every such map is a derivation of $\F_\epsilon^{\,2}$, it is
sufficient to check the relations (\ref{C-D}) for $G_S$ where $S$ is one
of the elementary traceless symmetric maps of the form $S(a,b)$. This is
straightforward. 

It is also straightforward (though tedious) to check that the maps
$E_a,F_a$ are derivations of $\F_\epsilon^{\,2}$ for any $a\in\F^{\,\prime}$. Now
let $D$ be any derivation of $\F_\epsilon^{\,2}$, and write 
\[
 Dl = \alpha + a +l(\beta + b)
\]
with $\alpha, \beta \in \R$ and $a,b\in\F^{\,\prime}$. Since $l^2=\epsilon$ and
$D\epsilon = 0$, $Dl$ must anticommute with $l$; hence $\alpha = \beta =
0$, so that 
\[
 Dl = -\frac{1}{2\epsilon}F_a(l) + E_b(l).
\]
It follows, by the first part of the proof, that $D +
(2\epsilon)^{-1}F_a - E_b$ is the sum of a derivation of $\F$ and an
element $G_S$. Thus the derivations $\overline{D}, E_a, F_a$ and $G_S$
span $\Der\F_\epsilon^{\,2}$.

The stated commutators can be verified by straightforward computation.
\end{proof}

Let $\Der_0\F_\epsilon^{\,2}$ be the subalgebra of $\Der\F_\epsilon^{\,2}$
spanned by $\Der\F$ and $E_a$ ($a\in\F^{\,\prime}$), so that 
\be
 D\in\Der_0\F_\epsilon^2 \iff D(\F)\subset \F \text{ and }
D(l\F)\subset l\F ; \label{Der-c}
\end{equation}
and let $\Der_1\F_\epsilon^{\,2}$ be the subspace spanned by $F_a$ and
$G_S$, so that 
\be
 D\in\Der_1\F_\epsilon^{\,2} \iff D(\F)\subset l\F \text{ and }
D(l\F)\subset\F. \label{Der-nc}
\end{equation}
Then $\Der\F_\epsilon^{\,2}$ has
the structure \eqref{comp-decomp}, with $\n = \Der_0\F_\epsilon^{\,2}$ and
$\p = \Der_1\F_\epsilon^{\,2}$, and the brackets in $[\p,\p]$, which include
a factor $\epsilon$, have opposite signs in $\K$ and $\sK$. Since
$\Der\K$ is compact (being a subalgebra of $\so(\K)$), this identifies
the maximal compact subalgebra of $\Der\sK$ as $\Der_0\sK =
\Der\F\dsum\F^{\,\prime}$; explicitly,
\[
 \Der_0\sC = 0, \qquad \Der_0\sH = \so(2), \qquad \Der_0\sO = \so(4).
\]

Note that the algebra $\F_\epsilon$ has a $\Z_2$-grading $\F_\epsilon^{\,2}
= \F\dsum l\F$, and the above decomposition is the corresponding
$\Z_2$-grading of the derivation algebra, i.e.\ $\Der_\delta\F_\epsilon^{\,2}$
($\delta = 0,1$) is the subspace of derivations of degree $\delta$. From
the definition \eqref{D-def} it follows that the derivation $D_{x,y}$
has degree $\gamma + \delta$~(mod~2) if $x$ has degree $\gamma$ and $y$
has degree $\delta$.

Now consider the rows of the non-compact magic square $L_3(\sK_1,\K_2)$.
Suppose $\sK_1 = (\F_2)_+^2 = \F_2\dsum l\F_2$. Then Vinberg's
construction gives
\begin{align*}
 L_3(\sK_1,\K_2) &= A_3\pr(\sK_1\ox \K_2)\dsum\Der\sK_1\dsum\Der\K_2\\
&=\n\dsum\p
\end{align*}
where
\begin{align*}
\n &= A_3\pr(\F_1\ox\K_2)\dsum\Der_0\sK_1\dsum\Der \K_2,\\
\p &= A_3(l\F_1\ox \K_2)\dsum\Der_1\sK_1.
\end{align*}
The brackets (\ref{Vin1}--\ref{Vin2}), together with the $\Z_2$-grading
of $\sK_1$ and $\Der\sK_1$, give the structure \eqref{comp-decomp}. The
compact algebra $L_3(\K_1\ox\K_2)$ has the same structure, and the
brackets are the same except for the sign in $[\p,\p]$, which contains a
factor $\epsilon$. Hence the maximal compact subalgebra of
$L_3(\sK_1,\K_2)$ is 
\begin{align*} \n &=
A_3\pr(\F_1\ox\K_2)\dsum\Der\K_2\dsum\Der\F_1\dsum\F_1^{\,\prime}\\
&=L_3(\F _1,K_2)\dsum\F_1^{\,\prime}.
\end{align*}
Thus we have

\begin{theorem} The maximal compact subalgebra of the non-compact magic
square algebra $L_3(\sK_1,K_2)$ is $L_3(\F_1,K_2)\dsum \F_1\pr$, where $\F_1$
is the division algebra preceding $\K_1$ in the Cayley-Dickson process.
\end{theorem}

Applying this to the last row and column of the magic square gives the table
at the end of Section \ref{sec:Tits}.

For completeness, we identify the real Lie algebras occurring in the
magic square $L_3(\sK_1,\sK_2)$ when both composition algebras are
split. Writing $\sK_i = (\F_i)_+^2$ ($i$=1,2) gives a
$(\Z_2\times\Z_2)$-grading
\begin{align*} L_3(\sK_1,\sK_2) = &A_3\pr(\F_1\ox\F_2) + \Der_0\sK_1
+\Der_0\sK_2\\ 
&+ A_3\pr(l\F_1\ox\F_2)\dsum\Der_1\sK_1\\
&+ A_3\pr(\F_1\ox l\F_2)\dsum\Der_1\sK_2\\
&+A_3\pr(l\F_1\ox l\F_2)
\end{align*}
in which the successive lines have gradings (0,0), (1,0), (0,1) and
(1,1). By arguments similar to those used for 
$L_3(\sK_1,\K_2)$, the maximal compact subalgebra is the direct sum of
the subspaces of degree (0,0) and (1,1), namely
\begin{align*}
\n &= A_3\pr(\F_1\ox\F_2)\dsum\Der\F_1 \dsum\F_1^{\,\prime} \dsum \Der\F_2 \dsum
\F_2^{\,\prime} \dsum A_3\pr(l_1\F_1\ox l_2\F_2)\\
&= L_3(\F_1,\F_2) \dsum \F_1^{\,\prime} \dsum \F_2^{\,\prime} \dsum A_3\pr(l_1\F_1\ox l_2\F_2).
\end{align*}
Since the elements of $l_1\F_1\ox l_2\F_2$ are self-conjugate in
$\K_1\ox\K_2$, the last summand contains antisymmetric $3\times 3$
matrices which can be identified with the entries in the last row and column
(excluding the diagonal element) of an antihermitian $4\times 4$ matrix
over $\F_1\ox\F_2$, while an element of $\F_1^{\,\prime}\dsum\F_2^{\,\prime}$ can be
identified with the last diagonal element of such a matrix. Thus we have
a vector space isomorphism
\be \label{4-square}
\n = L_4(\F_1,\F_2).
\end{equation}
We will find that this is actually a Lie algebra isomorphism.

By inspection of the table of real forms of complex semi-simple Lie
algebras \cite{Gilmore, Onishchik} we can now identify the non-compact
Lie algebras of the doubly-split magic square $L_3(\sK_1,\sK_2)$ as
follows:
\begin{center}
\vspace{0.5cm} \begin{tabular}{|c||c|c|c|c|}
\hline
  & $\mr$ & $\mc$  & $\mh$ & $\mo$ \\
 \hline \hline
  $\mr$ & $\so(3)$ & $\sl(3,\R)$ & $\sp(6,\R)$ & $F_4(-4)$ \\
  \hline
   $\mc$ & $\sl(3,\R)$ & $\sl(3,\R) \oplus \sl(3,\R)$ & $\sl(6,\R)$ & $E_6(-6)$ \\
   \hline
    $\mh$ & $\sp(6,\R)$ & $\sl(6,\R)$ & $\so(6,6)$ & $E_7(-7)$ \\
    \hline
    $\mo$ & $F_4(-4)$ & $E_6(-6)$ & $E_7(-7)$ & $E_8(-8)$ \\
    \hline
    \end{tabular}. \vspace{0.5cm}
    \end{center}

\noindent in which the real forms of the exceptional Lie algebras are identified
by the signatures of their Killing forms. The maximal compact subalgebras are 

\begin{center}
\vspace{0.5cm} \begin{tabular}{|c||c|c|c|c|}
\hline
  & $\mr$ & $\mc$  & $\mh$ & $\mo$ \\
 \hline \hline
  $\mr$ & $\so(3)$ & $\so(3)$ & $\u(3)$ & $\sq(3)\osum\so(3)$ \\
  \hline
   $\mc$ & $\so(3)$ & $\so(3) \oplus \so(3)$ & $\so(6)$ & $\sq(4)$ \\
   \hline
    $\mh$ & $\u(3)$ & $\so(6)$ & $\so(6)\osum\so(6)$ & $\su(8)$ \\
    \hline
    $\mo$ & $\sq(3)\osum\so(3)$ & $\sq(4)$ & $\su(8)$ & $\so(16)$ \\
    \hline
    \end{tabular}. \vspace{0.5cm}
    \end{center}

In this last table the $3\times 3$ square labelled by $\C$, $\H$ and
$\OO$ is isomorphic to 

\begin{center}
\vspace{0.5cm} \begin{tabular}{|c||c|c|c|}
\hline
   & $\mc$  & $\mh$ & $\mo$ \\
 \hline \hline
   $\mc$ & $\so(4)$ & $\su(4)$ & $\sq(4)$ \\
   \hline
    $\mh$ & $\su(4)$ & $\su(4)\osum\su(4)$ & $\su(8)$ \\
    \hline
    $\mo$ & $\sq(4)$ & $\su(8)$ & $\so(16)$ \\
    \hline
    \end{tabular}. \vspace{0.5cm}
    \end{center}
 
\noindent which has a non-compact form

\begin{center}
\vspace{0.5cm} \begin{tabular}{|c||c|c|c|}
\hline
   & $\mr$  & $\mc$ & $\mh$ \\
 \hline \hline
   $\mr$ & $\so(4)$ & $\su(4)$ & $\sq(4)$ \\
   \hline
    $\mc$ & $\sl(4,\R)$ & $\sl(4,\C)$ & $\sl(4,\H)$ \\
    \hline
    $\mh$ & $\sp(8,\R)$ & $\sp(8,\C)$ & $\sp(8,\H)$ \\
     &&$\cong\su(4,4)$&\\
    \hline
    \end{tabular}. \vspace{0.5cm}
    \end{center}

\noindent in which we have changed the labels of the rows and columns from $\K$ to
$\F$ where $\K =\F_\epsilon^{\,2}$ with $\epsilon = +1$ for the rows and
$\epsilon = -1$ for the columns. The rows of this table are
$\sa(4,\F_2)$, $\sl(4,\F_2)$ and $\sp(8,\F_2)$, and therefore by theorem
\ref{thm:ntimesn} the Lie algebras in the table are $L_4(\sF_1,
\F_2)$. The compact
forms are therefore $L_4(\F_1,\F_2)$ as asserted in \eqref{4-square},
and we have established that this is a Lie algebra isomorphism.

The involution of the compact Lie algebra $L_3(\K_1,\K_2)$ which defines
the non-compact form $L_3(\sK_1,\sK_2)$ can be taken to be $X\mapsto
-X^T$ for $X\in A_3\pr(\sK_1,\sK_2)$, together with the (essentially
unique) non-trivial involution on both $\Der\K_1$ and $\Der\K_2$. The
Cartan subalgebra of $L_3(\K_1,\K_2)$ can be chosen so that this
involution takes each root element $x_\alpha$ to $x_{-\alpha}$ (and
preserves the Cartan subalgebra). This explains why the rank of each of
the Lie algebras $L_3(\sK_1,\sK_2)$ is equal in magnitude to the signature of its
Killing form.

The magic squares $L_3(\K_1,\sK_2)$ and $L_3(\sK_1,\sK_2)$ contain all
the real forms of the exceptional simple Lie algebras except the
following two:
\begin{align*} F_4(20) &\text{ with maximal compact subalgebra }
\so(9);\\ 
E_6(14) &\text{ with maximal compact subalgebra } \so(10)\oplus\so(2).
\end{align*}
These can presumably be explained by a construction in which the antihermitian
matrices $A_3\pr(\K_1\otimes\K_2)$ are replaced by matrices which are
antihermitian with respect to a non-positive definite metric matrix $G =
\text{diag}(1,1,-1)$, i.e.\ by matrices $X$ satisfying
$\overline{X}^TG=-GX$, as in the Santander-Herranz construction.

\section{The $n = 2$ Magic Square}
\label{2mat}

It would be surprising, particularly in view of Freudenthal's
geometrical interpretation (see Section \ref{rows}), if $n=3$ were the 
only case in which there were
Lie algebras $L_n(\K_1,\K_2)$ for non-associative $\K_1$ and $\K_2$; we
would expect Lie algebras corresponding to $n=2$ to arise as subalgebras
of the $n=3$ algebras. Indeed, the algebra of $2\times 2$
hermitian matrices $H_2(\K)$ is a Jordan algebra if $H_3(\K)$ is, and
therefore the Tits construction of Theorem~\ref{thm:Tits} yields a Lie
algebra $L_2(\K_1,\K_2)=L(\K_1,H_2(\K_2))$ for associative $\K_1$ and
for any composition algebra $\K_2$. We will now show how to extend this
construction to allow $\K_1$ to be any composition algebra.

For $n=2$ the hermitian Jordan algebra $H_2(\K)$ takes a particularly
simple form. The usual identification of $\R$ with the subspace of scalar
multiples of the identity gives $H_2(\K) = \R \dsum H_2\pr(\K)$, and the
Jordan product in the traceless subspace $H_2\pr(\K)$ is given by 
\be
 A\cdot B = \<A,B\>\1 \label{JorCliff}
\end{equation}
where the inner product is defined by 
\[
 \< A,B\> = \half\tr (A\cdot B) 
\]
so that 
\be
 A = \begin{pmatrix}\lambda & x\\\overline{x}&-\lambda\end{pmatrix},
 B = \begin{pmatrix}\mu & y\\\overline{y}&-\mu\end{pmatrix}
\implies \<A,B\> = 2\big(\lambda\mu + \<x,y\>\big) \label{H2'innerp},
\end{equation}
i.e.
\be
 H_2\pr(\K) = \R\osum\K \label{H2'}
\end{equation}
(recall that we use $\osum$ to denote that the summands are orthogonal
subspaces). The anticommutator algebra of $H_2(\K)$ is therefore a
subalgebra of that of the Clifford algebra of the vector space
$\R\osum\K$; since the Clifford algebra is associative, its
anticommutator algebra is a (special) Jordan algebra. It is immediate
from \eqref{JorCliff} that the derivations of this Jordan algebra are
precisely the antisymmetric linear endomorphisms of $H_2\pr(\K)$. To
summarise, 

\begin{theorem}\label{H2} If $\K$ is any composition algebra, the
anticommutator algebra $H_2(\K)$ is a Jordan algebra with product given
by \eqref{JorCliff}, and its derivation algebra is 
\be
 \Der H_2(\K) \cong \so(\R\osum\K). \label{derH2}
\end{equation}
\end{theorem}

There is also a description of $\Der H_2(\K)$ in terms of $2\times 2$
matrices like, but interestingly different from, the description of
$\Der H_3(\K)$ in Theorem~\ref{DerH3}:

\begin{theorem}\label{DerH2} For any composition algebra $\K$,
\[
 \Der H_2(\K) = \so(\K\pr)\dsum A_2\pr(\K)
\]
in which $\so(\K\pr)$ is a Lie subalgebra, the Lie brackets between 
$\so(\K\pr)$ and $A_2^{\prime}(\mk)$  are given by the elementwise action of
$\so(\K\pr)$ on $2\times 2$ matrices over $\K$, and
\begin{equation}\label{A2brac}
    [A,B]=(AB-BA)^{\prime}+2S(A,B)
\end{equation}
where $A,B \in A_2^{\prime}(\mk)$, the prime denotes the traceless
part, and
\begin{equation*}
    S(A,B)=\sum_{ij} S_{a_{ij},b_{ij}}\in \so(\K\pr).
\end{equation*}
\end{theorem}

\begin{proof}
Write $H_2\pr(\K) = \sigma_1(\K\pr)\dsum\R\sigma_2\dsum\R\sigma_3$ where
$\sigma_1:\K\pr\rightarrow H_2\pr(\K)$ and $\sigma_2,\sigma_3\in H_2\pr(\K)$
are defined by 
\[ \sigma_1(a) = \begin{pmatrix}0&a\\-a&0\end{pmatrix},
\quad \sigma_2 = \begin{pmatrix}0&1\\1&0\end{pmatrix},
\quad \sigma_3 = \begin{pmatrix}1&0\\0&-1\end{pmatrix}.
\]
Then 
\[
 \Der H_2(\K) = \so(H_2\pr(\K)) = \so(\K\pr) \dsum\R \theta_1 \dsum
\theta_2(\K\pr) \dsum \theta_3(\K\pr)
\]
where the actions of the derivations $\theta_1$, $\theta_2(a)$ and
$\theta_3(b)$ ($a,b\in\K\pr$) are given by
\begin{align*}
&\theta_1(a\sigma_1) = 0, &&\theta_1(\sigma_2) = \sigma_3,
&&\theta_1(\sigma_3) = -\sigma_2,\\
&\theta_2(a)(b\sigma_1) = -\<a,b\>\sigma_3, &&\theta_2(a)\sigma_2 = 0,
&&\theta_2(a)\sigma_3 = a\sigma_1,\\
&\theta_3(a)(b\sigma_1) = \<a,b\>\sigma_2, &&\theta_3(a)\sigma_2 =
-a\sigma_1, &&\theta_3(a)\sigma_3 = 0.
\end{align*}
These actions are reproduced by
\begin{align} \theta_1(H) &= [\sigma_1, H] \notag\\
\theta_2(a)(H) &= [a\sigma_2, H] \label{A2com}\\
\theta_3(a)(H) &= [a\sigma_3, H] \notag
\end{align}
where the brackets denote matrix commutators, so
\[
 \R\theta_1 \dsum \theta_2(\K\pr) \dsum \theta_3(\K\pr) = 
\R\sigma_1 \dsum \K\pr\sigma_2 \dsum \K\pr\sigma_3 = A_3\pr(\K),
\]
and the action of $A\in A_3(\K)$ as a derivation of $H_2(\K)$ is the
commutator map $C_A$. By the matrix identity \eqref{CACB2}, the Lie
bracket in $\Der H_2(\K)$ is given by
\[
    [A,B]=(AB-BA)^{\prime}+\half F(A,B).
\]
where $F(A,B)$ is defined in \eqref{F-def}.
Now for $A,B\in A_2(\H)$ and $z\in\K\pr$ we have
\begin{align*}
4S(A,B)z &= \sum_{ij}\big\{ \<a_{ij},z\>b_{ji} - \<b_{ji},z\>a_{ij}\big\}\\
&= \sum_{ij}\big\{ (a_{ij}\overline{z} + z\overline{a_{ij}})b_{ji} +
b_{ji}(\overline{a_{ij}}z + \overline{z}a_{ij})\\ 
&\qqqquad - (b_{ji}\overline{z} -
z\overline{b_{ji}})a_{ij} - a_{ij}(\overline{b_{ji}}z +
\overline{z}b_{ji})\big\}\\
&= \sum_{ij}\big\{ z(\overline{a_{ij}}b_{ji} -
\overline{b_{ji}}a_{ij}) - (b_{ji}\overline{a_{ij}} -
a_{ij}\overline{b_{ji}})z\big\}\\ 
&= F(A,B)z.
\end{align*}
The Lie bracket can therefore be written as \eqref{A2brac}.
\end{proof}

Comparison with Theorem~\ref{DerH3} suggests that in passing from $n=3$
to $n=2$, $\Der\K$ should be replaced by $\so(\K\pr)$. This has no
effect if $\K$ is associative, since then these two Lie algebras
coincide\footnote{This seems to be a genuine coincidence since it does
not survive at the group level: Aut$\C$ = O($\C\pr$) but Aut$\H$ =
SO($\H\pr$).}  (see Section~\ref{notation}). Thus we will define
$L_2(\K_1,\K_2)$ by making this replacement in the definition of
$L_n(\K_1,\K_2)$ for $n\ge 3$, which gives the vector space
\be \label{L2def}
 L_2(\K_1,\K_2) = \so(\K_1\pr) \dsum \Der H_2(\K_2) \dsum
\K_1\pr\otimes H_2\pr(\K_2) 
\end{equation}
To obtain the replacement for the Lie bracket \eqref{Titsbrac}, note that
it follows from \eqref{JorCliff} that the traceless Jordan product $A*B$
is identically zero in $H_2(\K)$. Moreover, if $\K_1$ is associative the
derivation $D_{a,b}\in\Der\K_1$ is the generator $4S_{a,b}$
of rotations in the plane of $a$ and $b$, for 
\begin{align*}
4S_{a,b}(c) &= 2\<a,c\>b + 2b\<a,c\> -2\<b,c\>a - 2a\<b,c\>\\
&= -(ac+ca)b - b(ac+ca) + (bc+cb)a + a(bc+cb)\\ 
&=[[a,b],c].
\end{align*}
Thus the following definition of $L_2(\K_1,\K_2)$ is the case $n=2$ of
the general definition of $L_n$ if $\K_1$ is associative:

\begin{defn} The algebra $L_2(\K_1,\K_2)$ consists of the vector space
\eqref{L2def} with brackets in the first two summands given by the Lie
algebra $\so(\K_1)\oplus\Der H_2(\K_2)$, brackets between these and the
third summand given by the usual action on $\K_1\pr\otimes
H_2\pr(\K_2)$, and further brackets
\be
[a\otimes A,b\otimes B] = 2\<A,B\>S_{a,b}  - 4\<a,b\>[L_A,L_B]
\label{L2brac} 
\end{equation}
($A,B\in H_2\pr(\K_2)$).
\end{defn}

This is readily identified as a Lie algebra:

\begin{theorem}\label{L2} If $\K_1$ and $\K_2$ are composition
algebras, $L_2(\K_1,\K_2)$ as defined above is a Lie algebra isomorphic
to $\so(\K_1\otimes\K_2)$.
\end{theorem}

\begin{proof} The orthogonal Lie algebra $\so(V)$ of a real
inner-product vector space $V$ is spanned by the elementary generators
$S_{x,y}$, defined as in \eqref{S-def} with $x,y\in V$. Hence 
\[
\so(V\oplus W) \cong \so(V) \dsum \so(W) \dsum S_{V,W}
\]
where $S_{V,W}$, spanned by $S_{v,w}=-S_{w,v}$ with $v\in V$, $w\in W$,
is isomorphic to $V\otimes W$.  Taking $V=\K_1\pr$,
$W=\R\osum\K_2\cong H_2\pr(\K_2)$, we have a vector space isomorphism $\theta$
between $\so(V\osum W) \cong \so(\K_1\osum\K_2)$ and $L_2(\K_1,\K_2)$
such that $\theta|\so(\K_1\pr)$ is the identity,
$\theta|\so(\R\osum\K_2)$ is the isomorphism of Theorem~\ref{H2}, and
$\theta|S_{V,W}$ is given by
\[
\theta(S_{a,A}) = \frac{1}{\sqrt{2}}a\otimes A \qquad 
(a\in\K_1\pr,\; A\in H_2\pr(\K_2\pr)).
\]
Then $\theta$ is an algebra isomorphism: for $X\in\so(\K_1\pr)$,
\[
\theta([X,S_{a,A}]) = \theta(S_{Xa,A}) = 2(Xa\otimes A) =
[\theta(X),\theta(S{a,A})], 
\]
and similarly for $Y\in\so(\R\osum\K_2)$. Finally,
\[
\theta\big([S_{a,A},S_{b,B}]) = \theta(\<a,b\>S_{A,B} + \<A,B\>S_{a,b}\big)
\]
while
\[
[\theta(S_{a,A}),\theta(S_{b,B})] = \<A,B\>S_{a,b} - 2\<a,b\>[L_A,L_B]
\]
and \eqref{JorCliff} gives $[L_A,L_B] = -\half S_{A,B}$. Hence
$L_2(\K_1,\K_2)$ is a Lie algebra and $\theta$ is a Lie algebra
isomorphism.
\end{proof}

This theorem shows that the compact and doubly split magic squares
$L_2(\K_1,\K_2)$ and $L_2(\sK_1,\sK_2)$ ($\K_1,\K_2=\R,\C,\H,\OO$), like
the $n=3$ squares, are symmetric. The complex types of these Lie
algebras are identified below.
\begin{center}
{\renewcommand{\arraystretch}{1.2}
\vspace{0.5cm} \begin{tabular}{|c||c|c|c|c|} \hline
    & $\mr$ & $\mc$ & $\mh$ & $\mo$ \\ \hline \hline
$L_2(\mk ,\mr )$ & $D_{1}$ & $A_{1}\cong B_{1} \cong C_{1}$ &
$C_{2}\cong B_{2}$ & $B_{4}$ \\ \hline
$L_2(\mk ,\mc )$ & $A_{1}\cong B_{1} \cong C_{1}$ & $A_{1} \oplus
A_{1}$ & $A_{3}\cong D_{3}$ & $D_{5}$ \\ \hline
$L_2(\mk ,\mh )$ & $C_{2}\cong B_{2}$ & $A_{3}\cong D_{3}$ & $D_{4}$
& $D_{6}$ \\ \hline
$L_2(\mk ,\mo )$ & $B_{4}$ & $D_{5}$ & $D_{6}$ & $D_{8}$ \\ \hline
\end{tabular}} \vspace{0.5cm}.
\end{center}
The following table shows the mixed square $L_2(\sK_1,\K_2)$. 
\begin{center}
{\renewcommand{\arraystretch}{1.4}
\vspace{0.5cm} \begin{tabular}{|c||c|c|c|c|} \hline
    & $\mr$ & $\mc$ & $\mh$ & $\mo$ \\ \hline \hline
$L_2(\mk ,\mr )$ & $\so (2)$ & $\so (3)$ & $\so (5)$ & $\so (9)$ \\ 
\hline\vspace{2pt}
$L_2(\mcs ,\mk )$ & $\so (2,1)$ & $\so (3,1)$ & $\so (5,1)$ & $\so (9,1)$ \\
\hline\vspace{2pt}
$L_2(\mhs ,\mk )$ & $\so (3,2)$ & $\so (4,2)$ & $\so (6,2)$ & $\so (10,2)$ \\
\hline\vspace{2pt}
$L_2(\mos ,\mk )$ & $\so (5,4)$ & $\so (6,4)$ & $\so (8,4)$ & $\so (12,4)$ \\
\hline
\end{tabular}} \vspace{0.5cm}.
\end{center}
Note that the maximal compact subalgebras in each row of this square are
related to the previous row as in the $n=3$ magic square
(Theorem~\ref{maxcomp}). Because the definition of $L_2(\K_1,\K_2)$
coincides with the Tits construction $T(\K_1,H_2(\K_2))$ if $\K_1$ is
associative, Theorem~\ref{Jmatrix} gives the same relation between
the rows of the non-compact $n=2$ square $L_2(\sK_1,\K_2)$ and the
matrix Lie algebras $\sa(2,\K)$, $\sl(2,\K)$ and $\sp(4,\K)$ as for
$n=3$ as shown below.
\begin{center}
{\renewcommand{\arraystretch}{1.4}
\vspace{0.5cm} \begin{tabular}{|c||c|c|c|c|} \hline
    & $\mr$ & $\mc$ & $\mh$ & $\mos$ \\ \hline \hline \vspace{1pt}
$\Der H_2(\mk )\cong L_2(\mr ,\mk )$ & $\so(2)$ & $\su(2)$ & $\sq(2)$ 
& $\so(9)$ \\ \hline\vspace{2pt}
$\Str  H_2(\mk )\cong L_2(\mcs ,\mk )$ & $\sl(2,\mr )$ & $\sl(2,\mc )$ &
$\sl(2,\mh )$ & $\sl(2,\mo )$ \\ \hline\vspace{2pt}
$\Con  H_2(\mk )\cong L_2(\mhs ,\mk )$ & $\sp(4,\mr )$ & $\su(2,2)$ &
$\sp(4,\mh )$ & $\sp(4,\mo )$ \\ \hline\vspace{2pt}
$L_2(\mos ,\mk )$ & $\so(5,4)$ & $\so(6,4)$ & $\so(8,4)$ & $\so(12,4)$ \\ \hline
\end{tabular}} \vspace{0.5cm}.
\end{center}
The differences between the definitions of $L_2$ and $L_3$, however,
affect the description of these matrix Lie algebras as follows:

\newpage
\begin{theorem}\label{2matrix}\cite{Sudbery84}
For any composition algebra $\K$,

\linespace
\noindent \emph{\textbf{(a)}}\hspace{3cm}
$\sl(2,\K) = A_2\pr(\K)\dsum\so(\K\pr)$;

\linespace
\noindent\emph{\textbf{(b)}}\hspace{3cm} 
$\sl(2,\K) = M_2\pr(\K)\dsum\so(\K\pr)$;

\linespace
\noindent\emph{\textbf{(c)}}\hspace{3cm} 
$\sp(4,\K) = Q_4\pr(\K)\dsum\so(\K\pr)$.

\linespace

In each case  the Lie brackets are defined as follows:
\begin{enumerate}
\item $\so(\K\pr)$ is a Lie subalgebra;
\item The brackets between $\so(\K\pr)$ and the other summand are given by
the elementwise action of $\so(\K\pr)$ on matrices over $\K$;
\item The brackets between two matrices in the first summand are
\be \label{mat2brac}
[X,Y] = (XY-YX)\pr + \frac{1}{n}F(X,Y) 
\end{equation}
where $n$ (= 2 or 4) is the size of the matrix and $F(X,Y)\in\so(\K\pr)$ is defined
by
\[
F(X,Y)a = \sum_{ij}\big([[x_{ij},y_{ji}],a] + 2[x_{ij},y_{ji},a]\big)
\qquad (a\in\K\pr).
\]
\end{enumerate}
\end{theorem}
\begin{proof}
\textbf{(a)} The first isomorphism comes from Theorem~\ref{DerH2},
while the brackets are given by \eqref{A2brac}. 
This establishes \eqref{mat2brac} for $\sa(2,\K)$.

\linespace
\textbf{(b)} For $\sl(2,\K)$ the vector space is 
\begin{align*}
\sl(2,\K) = \text{Str}\pr H_2(\K) &= \Der H_2(\K) \dsum H_2\pr(\K)\\
&= \so(\K\pr)\dsum A_2\pr(\K) \dsum H_2\pr(\K)\\
&= \so(\K\pr)\dsum L_2\pr(\K).
\end{align*}
For $A,B\in H_2\pr(\K)$ the Lie bracket is that of $\Der H_2(\K)$, which we
have just seen to be given by \eqref{mat2brac}. For $A\in A_2\pr(\K)$,
$H\in H_2\pr(\K)$, the bracket is given by the action of $A$ as an element
of $\Der H_2(\K)$ on $H$, which according to \eqref{A2com} is
\be
[A,H] = AH - HA. \label{AHcom}
\end{equation}
Now by Lemma~\ref{trcomAH}, $\tr(XH-HX)=0$ and $F(A,H)=0$. 
Hence \eqref{AHcom} is the same as \eqref{mat2brac}.

Finally, for $H,K\in H_2\pr(\K)$ the $\Str\pr H_3(\K)$ bracket is 
\[
[H,K] = L_H L_K - L_K L_H \in \Der H_2(\K)
\]
and by Lemma~\ref{LHLK2} this commutator appears in the 
decomposition $\Der H_2(\K)=A_2\pr(\K)\dsum\so(\K\pr)$ as
\[ [L_H,L_K] = (HK-KH)\pr + \half F(H,K).
\]

As in Theorem~\ref{3x3mat}, the action of the matrix part of 
Str$\pr H_2(\K)$ on $H_2(\K)$ is
\be \label{StrH2act}
H \mapsto XH + HX^\dagger \qqquad (X\in M_2\pr (\K)).
\end{equation}

\linespace
\textbf{(c)} The proof of (c) is the same as that of
Theorem~\ref{3x3mat}(c) with the derivation $D(X,Y)$ replaced by the
orthogonal map $F(X,Y)$.
\end{proof}

\appendix
\section{Matrix identities}\label{appendix}

In this appendix we prove various identities for matrices with entries
in a composition algebra $\K$. For associative algebras these are
familiar Jacobi-like identities and trace identities; in general they
hold only for certain classes of matrix, and some need to be modified by
terms containing the following elements of $\Der\K$ and $\so(\K\pr)$
defined for pairs of $3\times 3$ matrices $X,Y$:
\be\label{Dmatdef}
D(X,Y) = \sum_{ij}D_{x_{ij},y_{ji}}
\end{equation}
where $D_{x,y}$ is the derivation of $\K$ defined in \eqref{D-def};
\be\label{Smatdef}
S(X,Y) = \sum_{ij}S_{x_{ij},y_{ji}}
\end{equation}
where $S_{x,y}$ is the generator of rotations in the plane of $x$ and
$y$, defined in \eqref{S-def};
\begin{equation}
E(X,Y)z = \sum_{ij}[x_{ij},y_{ji},z] \qquad (z\in\K),\label{E-def}
\end{equation}
and
\begin{equation} 
F(X,Y)z = \sum_{ij}[[x_{ij},y_{ji}],z] + 2[x_{ij},y_{ji},z].\label{F-def}
\end{equation}

In all of the following identities $\K$ is a composition algebra. The
square brackets denote matrix commutators and the chain brackets denote 
matrix anticommutators:
\[ 
[X,Y] = XY-YX, \qqquad \{X,Y\} = XY + YX.
\]

\begin{lemma}\label{3x3Jacobi}
The following matrix identities hold for $A,B\in A_3\pr(\K)$ and $H,K,L\in H_3(\K)$:

\linespace
\noindent\emph{\textbf{(a)}}\upline
\be\label{AHKJacobi}
[A,\{H,K\}] = \{[A,H],K\} + \{H,[A,K]\},
\end{equation}
\emph{\textbf{(b)}}\upline
\be\label{ABHJacobi}
[A,[B,H]]-[B,[A,H]] = [[A,B],H] + E(A,B)H,
\end{equation}
\emph{\textbf{(c)}}\upline
\be\label{HKLJacobi}
\{H,\{K,L\} - \{K,\{H,L\}\} = [[H,K],L] + E(H,K)L.
\end{equation}

\begin{proof}
\textbf{(a)} The difference between the two sides of \eqref{AHKJacobi} 
can be written in terms of matrix associators, whose $(i,j)$th element is
\begin{multline}
\label{eqn:paddington}
\sum _{mn}\big([a_{im},h_{mn},k_{nj}]+[a_{im},k_{mn},h_{nj}]+[h_{im},k_{mn},a_{nj}]\\
+[k_{im},h_{mn},a_{nj}]-[h_{im},a_{mn},k_{nj}]-[k_{im},a_{mn},h_{nj}]\big).
\end{multline}
Suppose $i\neq j$ and let $k$ be the third index. Since the diagonal
elements of $H$ and $K$ are real, any associator containing them
vanishes. Hence the terms containing $a_{ij}$ or $a_{ji}$ are
\begin{multline*}
    \sum_{n}\big([a_{ij},h_{jn},k_{nj}]+[a_{ij},k_{jn},h_{nj}]\big) + 
\sum_{m} \big([h_{im},k_{mi},a_{ij}]+ [k_{im},h_{mi},a_{ij}]\big) \\
-[h_{ij},a_{ji},k_{ij}]+[k_{ij},a_{ji},h_{ij}]=0
\end{multline*}
by the alternative law, the hermiticity of $H$ and $K$, and the fact
that an associator changes sign when one of its elements is conjugated.
The terms containing $a_{ik}$ or $a_{ki}$ are
\begin{equation*}
[a_{ik},h_{ki},k_{ij}]+[a_{ik},k_{ki},h_{ij}]-[h_{ik},a_{ki},k_{ij}]-
[k_{ik},a_{ki},h_{ij}]=0
\end{equation*}
using also $a_{ki}=-\overline{a}_{ik}$. Similarly, the terms containing
$a_{jk}$ or $a_{kj}$ vanish. Finally, the terms containing
$a_{ii},a_{jj}$ and $a_{kk}$ are
\begin{multline*}
[a_{ii},h_{ik},k_{kj}]+[a_{ii},k_{ik},h_{kj}]+[h_{ik},k_{kj},a_{jj}]+[k_{ik},h_{kj},a_{jj}]
 \\ -[h_{ik},a_{kk},k_{kj}]-[k_{ik},a_{kk},h_{kj}] =0
\end{multline*}
since $a_{ii}+a_{jj}+a_{kk}=0$.

Now consider the $(i,i)$th element. The last two terms of
equation~(\ref{eqn:paddington}) become
\begin{equation*}
-\sum_{mn} \left( [h_{im},a_{mn},k_{ni}]+[k_{in},a_{nm},h_{mi}]\right)
=0.
\end{equation*}
Let $j$ be one of the other two indices. The terms containing $a_{ij}$
or $a_{ji}$ are
\begin{equation*}
[a_{ij},h_{jk},k_{ki}]+[a_{ij},k_{jk},h_{ki}]+[h_{ik},k_{kj},a_{ji}]+[k_{ik},h_{kj},a_{ji}]=0,
\end{equation*}
where $k$ is the third index. There are no terms containing $a_{jk}$ or
$a_{kj}$. The terms containing $a_{ii},a_{jj}$ or $a_{kk}$ are
\begin{multline*}
\sum_{n} \left( [a_{ii},h_{in},k_{ni}]+[a_{ii},k_{in},h_{ni}]\right) \\
+\sum_{m} \left( [h_{im},k_{mi},a_{ii}]+[k_{im},h_{mi},a_{ii}]\right)=0.
\end{multline*}
Thus in all cases the expression~(\ref{eqn:paddington}) vanishes, proving (a).

\linespace
\textbf{(b), (c)} Similar arguments establish equations \eqref{ABHJacobi}
and \eqref{HKLJacobi}.
\end{proof}
\end{lemma}

\linespace
\begin{lemma}
\label{lemma:2x2Jacobi}
The identities of Lemma~\ref{3x3Jacobi} hold for $A,B\in A_2(\K)$ and 
$H,K,L\in H_2(\K)$.

\begin{proof}
The $2\times 2$ case can be deduced from Lemma~\ref{3x3Jacobi} by
applying it to the matrices  
\[
    \widetilde{A}=\begin{pmatrix} A & 0 \\ 0 & -\tr A \end{pmatrix},
\quad \widetilde{B}=\begin{pmatrix} B & 0 \\ 0 & -\tr B \end{pmatrix},
\]
\[
\quad \widetilde{X}=\begin{pmatrix} X & 0 \\ 0 & 1 \end{pmatrix},
\quad \widetilde{Y}=\begin{pmatrix} Y & 0 \\ 0 & 1 \end{pmatrix}.
\quad \widetilde{X}=\begin{pmatrix} Z & 0 \\ 0 & 1 \end{pmatrix}.
\]
Note that $A$ and $B$ do not need to be traceless.
\end{proof}
\end{lemma}

In the next identities we use the notation $L_H$ for the multiplication
by $H$ in the Jordan algebras $H_n(\K)$ ($n=2,3$) and $C_X$ for the
commutator with any matrix $X$:
\[ 
L_H(K) = \{H,K\}, \qquad C_X(Y) = [X,Y]
\]

\begin{lemma}\label{com3x3}
\emph{\textbf{(a)}} For $A,B\in A_3\pr(\K)$ and $H\in H_3(\K)$,
\be\label{CACB3}
[C_A, C_B]H = C_{(AB-BA)\pr}H+\third D(A,B)H.
\end{equation}
\emph{\textbf{(b)}} For $H,K,M\in H_3\pr(\K)$,
\be\label{LHLK3}
[L_H, L_K]M = C_{(HK-KH)\pr}M+\third D(H,K)M.
\end{equation}
\begin{proof}
By Lemma~\ref{3x3Jacobi},
\[ 
[C_A,C_B]H = C_{(AB-BA)\pr}H + [t_1,H] + E(X,Y)H
\]
and
\[
[L_H,L_K]M = C_{(HK-KH)\pr}M + [t_2,M] + E(H,K)M
\]
where 
\be\label{traces} 
t_1 = \third\tr(AB-BA) \quad \text{ and } \quad t_2 = \third\tr(HK-KH).
\end{equation}
But for any matrices $X,Y$,
\begin{align}\label{Dcalc}
[\third\tr(XY-YX),z] + E(X,Y)z &= \sum_{ij}\big(\third[[x_{ij},y_{ji}],z] 
+ [x_{ij},y_{ji},z]\big)\\ 
&= \third D(X,Y)z,
\end{align}
so the stated identities follow.
\end{proof}
\end{lemma}

\begin{lemma}\label{com2x2}
\emph{\textbf{(a)}} For $A,B\in A_2\pr(\K)$ and $H\in H_2(\K)$,
\be\label{CACB2}
[C_A, C_B]H = C_{(AB-BA)\pr}H+\half F(A,B)H.
\end{equation}
\emph{\textbf{(b)}} For $H,K,M\in H_2\pr(\K)$,
\be\label{LHLK2}
[L_H, L_K]M = C_{(HK-KH)\pr}M+\half F(H,K)M.
\end{equation}
\begin{proof}
The proof is the same as that of Lemma~\ref{com3x3} except that in
\eqref{traces} the fraction occurring is $\half$ rather than $\third$,
which means that in \eqref{Dcalc} $D(X,Y)$ must be replaced by $F(X,Y)$.
\end{proof} 
\end{lemma}

\begin{lemma}
\label{trcomAH}
For $A\in A_n(\K)$ and $H\in H_n(\K)$,
\[ 
\tr[A,H] = 0 \qquad \text{and} \qquad D(A,H) = F(A,H) = 0.
\]
\begin{proof}
\begin{align*}
\tr [A,H] = \sum_{ij}[a_{ij},h_{ji}] = \sum_{ij}
[\overline{a_{ij}}, \overline{h_{ji}}]\\
&= -\sum_{ij}[a_{ji},h_{ij}] = -\tr[A,H],
\end{align*}
so $\tr[A,H]=0$; and similarly,
\begin{align*}
D(A,H) &= \sum_{ij}D_{a_{ij},h_{ji}} =
\sum_{ij}D_{\overline{a_{ij}},\overline{h_{ji}}} =
\sum_{ij}D_{-a_{ji},h_{ij}}\\
&= - D(A,H),
\end{align*}
so $D(A,H)=0$. A similar argument shows that $F(A,H)=0$.
\end{proof}
\end{lemma}

\vspace{1cm}

\end{document}